%% file: paper.tex
\pdfoutput=1
\documentclass[letterpaper,twoside]{article}
\usepackage[preprint]{aistats2027}
\usepackage[hyphens]{url}
\usepackage{graphicx}
\usepackage[round,authoryear]{natbib}
\renewcommand{\bibsection}{\subsubsection*{\refname}}
\usepackage{caption}

\usepackage{amsmath}
\usepackage{amssymb}
\usepackage{amsthm}
\usepackage{algorithm}
\usepackage{algorithmic}
\usepackage{mathtools}
\usepackage{mathrsfs}
\usepackage{dsfont}
\usepackage{stmaryrd}
\usepackage{bm}
\usepackage{booktabs}
\usepackage{multirow}
\usepackage{array}
\usepackage{tabularx}
\usepackage{hhline}
\usepackage{subcaption}
\usepackage{enumitem}
\usepackage{upgreek}
\usepackage{xspace}
\usepackage{environ}
\usepackage{thm-restate}
\usepackage[dvipsnames]{xcolor}

\usepackage{tikz}
\usepackage{pgfplots}
\pgfplotsset{compat=newest}
\usepgfplotslibrary{groupplots,dateplot}
\usetikzlibrary{patterns,shapes.arrows,positioning,arrows,shapes,shapes.arrows,arrows.meta,trees,shapes.misc,shapes.geometric,decorations.pathreplacing,automata}
\tikzset{%
  >={Latex[width=1.5mm,length=2mm]},
  vertex/.style={draw,circle,inner sep=0mm,semithick,minimum width=4mm},
  point/.style = {circle, draw, inner sep=0.04cm,fill,node contents={}},
  uvertex/.style={outer sep=0},
  bidir/.style={<->,dashed, line width=0.25mm},
  dir/.style={->, line width=0.25mm},
  regime/.style={shape=rectangle,fill=black,inner sep=0pt,minimum size=3pt,draw},
  node distance=1cm,
  font=\scriptsize\sffamily
}
\pgfdeclarelayer{back}
\pgfsetlayers{back,main}
\makeatletter
\pgfkeys{%
  /tikz/on layer/.code={
    \def\tikz@path@do@at@end{\endpgfonlayer\endgroup\tikz@path@do@at@end}%
    \pgfonlayer{#1}\begingroup%
  }%
}
\makeatother

\theoremstyle{plain}
\newtheorem{theorem}{Theorem}

\newtheorem{lemma}[theorem]{Lemma}
\newtheorem{corollary}[theorem]{Corollary}
\theoremstyle{definition}

\theoremstyle{remark}

\theoremstyle{definition}
\newtheorem{example}{Example}

\definecolor{betterred}{RGB}{228,26,28}
\definecolor{betterblue}{RGB}{55,126,184}

\makeatletter
\newcommand{\xdashleftrightarrow}[2][]{\ext@arrow 3359\leftrightarrowfill@@{#1}{#2}}
\def\rightarrowfill@@{\arrowfill@@\relax\relbar\rightarrow}
\def\leftarrowfill@@{\arrowfill@@\leftarrow\relbar\relax}
\def\leftrightarrowfill@@{\arrowfill@@\leftarrow\relbar\rightarrow}
\def\arrowfill@@#1#2#3#4{%
  $\m@th\thickmuskip0mu\medmuskip\thickmuskip\thinmuskip\thickmuskip
   \relax#4#1
   \xleaders\hbox{$#4#2$}\hfill
   #3$%
}
\makeatother

\newcolumntype{C}[1]{>{\centering\arraybackslash}m{#1}}

\definecolor{plotred}{RGB}{163,74,90}  
\definecolor{plotdarkred}{RGB}{78,21,38}  
\definecolor{plotgray}{RGB}{115,115,115}  
\definecolor{plotblue}{RGB}{108,142,191}  
\definecolor{plotdarkblue}{RGB}{38,69,110}  
\definecolor{plotorange}{RGB}{217,160,91}  
\definecolor{plotdarkorange}{RGB}{199,123,37}  
\definecolor{plotgreen}{RGB}{111,158,130}  

\newcommand{\cblue}[1]{\textcolor{plotblue}{#1}}
\newcommand{\cdblue}[1]{\textcolor{plotdarkblue}{#1}}
\newcommand{\corange}[1]{\textcolor{plotorange}{#1}}
\newcommand{\cdorange}[1]{\textcolor{plotdarkorange}{#1}}

\newcommand{\Parens}[1]{\left(#1\right)}

\newcommand{\inv}[2]{P_{#2}\left ( #1\right)}
\newcommand{\invE}[2]{\3E_{#2}\left [ #1\right]}

\newcommand{\doo}{\text{do}}

\def\*#1{\boldsymbol{#1}}
\def\1#1{\mathcal{#1}}
\def\2#1{\mathscr{#1}}
\def\3#1{\mathbb{#1}}
\def\4#1{\mathds{#1}}
\def\5#1{\bar{\*#1}}

\newcommand{\papertitle}{Sharp Bounds on Continuous Treatment Effects under Unobserved Confounding}

\begin{document}
\runningtitle{Lipschitz Partial Identification of Conditional Treatment Effects}
\runningauthor{Zhang}

\twocolumn[
\aistatstitle{\papertitle}
\aistatsauthor{Junzhe Zhang}
\aistatsaddress{
Syracuse University \\
\texttt{jzhan403@syr.edu}
}
]

\input{section0}   
\input{section1}
\input{section2}
\input{section3}

\input{section4}
\input{section5}
\input{section6}

\label{sec:references_start}
\bibliography{references}

\appendix
\onecolumn
\thispagestyle{empty}
\aistatstitle{\papertitle: \\ Supplementary Materials}

\input{appendix_e}   
\input{appendix_a}
\input{appendix_b}

\end{document}

%% file: section0.tex
\begin{abstract}
We study partial identification of conditional treatment effects with a continuous treatment under unobserved confounding. We target the conditional average potential outcome (CAPO) $\mathbb{E}[Y_{\*x} \mid \*z]$ and the conditional average treatment effect (CATE) between two treatment levels. We assume a bounded outcome and Lipschitz continuity of the response functions and of the propensity scores given the latent confounder, and show that these conditions imply the smooth treatment response and smooth treatment selection conditions of Kim et al.\ (2018). The combined bound of Kim et al.\ uses only two reference points on the identified diagonal and is generally loose. Bounding through every diagonal reference point instead gives explicit bounds whose average is the sharp identified set of the CAPO at every treatment, inside or outside the observational support. Capping the difference between the upper and lower bounds gives sharp causal bounds on the CATE. These bounds are informative without monotonicity or overlap assumptions.
\end{abstract}

%% file: section1.tex
\section{Introduction}
Evaluating treatment effects from observational data is a central task in causal inference, with applications to policy evaluation, healthcare, and reinforcement learning. The challenge is that observed treatment--outcome correlations generally do not imply causal effects: unobserved confounding can induce spurious relations, leaving the target effect non-identifiable, i.e., underdetermined by the observational data together with causal knowledge about the underlying system \citep{pearl:2k}. Several conditions and algorithms have been developed to decide whether a causal effect is identifiable \citep{pearl:2k,bareinboim:pea:16-r450}. For example, under the \emph{back-door} condition \citep[Ch.~3.3.1]{pearl:2k}, or equivalently \emph{ignorability} \citep{rosenbaum:rub83}, the effect is identified by covariate adjustment and can be estimated at standard rates \citep{rosenbaum:rub83,bang2005doubly}. In practice, however, the available causal assumptions are often insufficient to determine the target effect from observational data, in which case the effect is said to be \emph{non-identifiable}. The following example illustrates this challenge. 

\begin{figure}[t]
\centering
    \hfill
        \begin{subfigure}{0.48\linewidth}\centering
            \setlength{\abovecaptionskip}{0pt}
		\includegraphics[width=\linewidth]{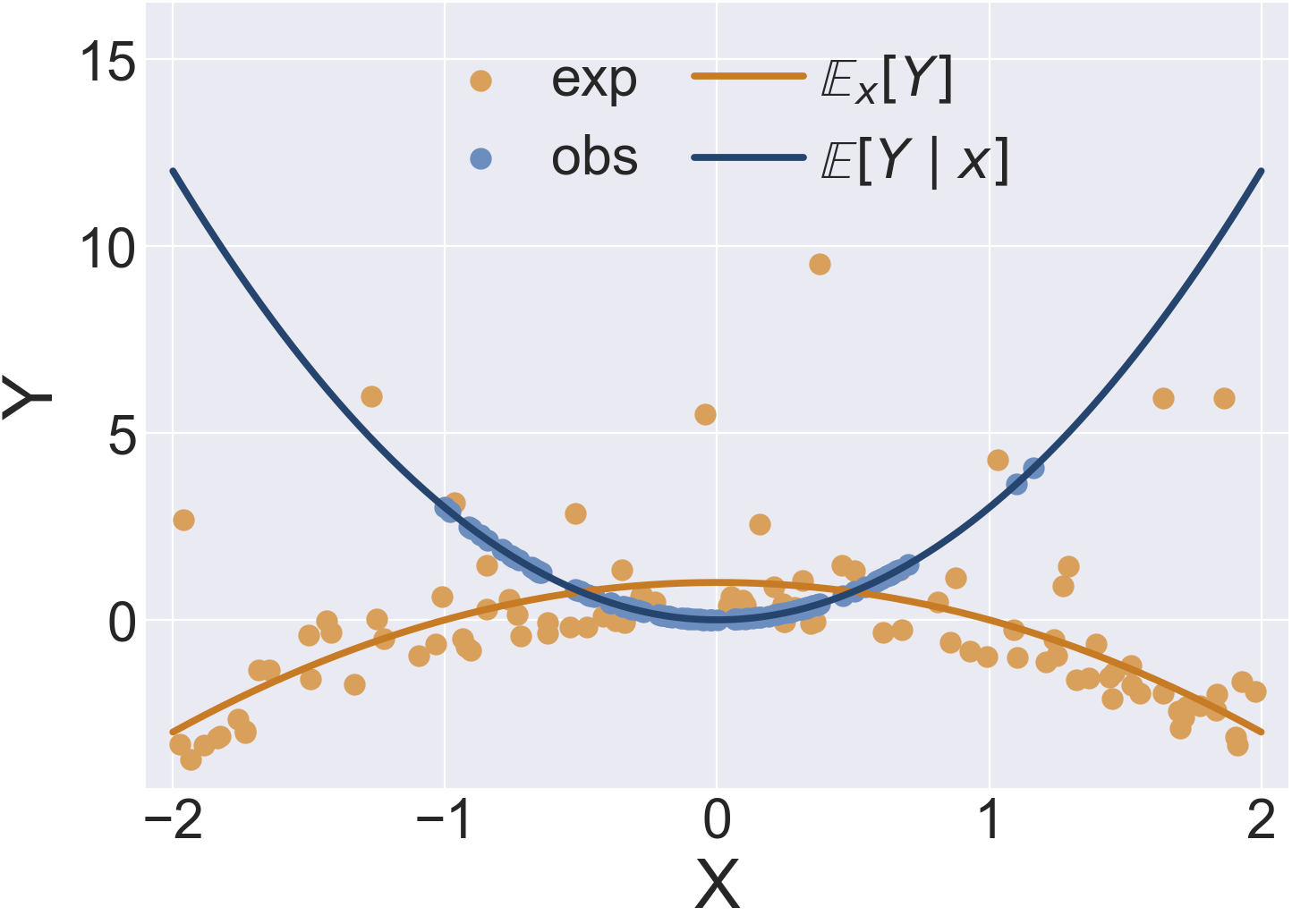}
		\caption{$\1M$}
		\label{fig:confounded}
	\end{subfigure}\hfill
        \begin{subfigure}{0.48\linewidth}\centering
            \setlength{\abovecaptionskip}{0pt}
		\includegraphics[width=\linewidth]{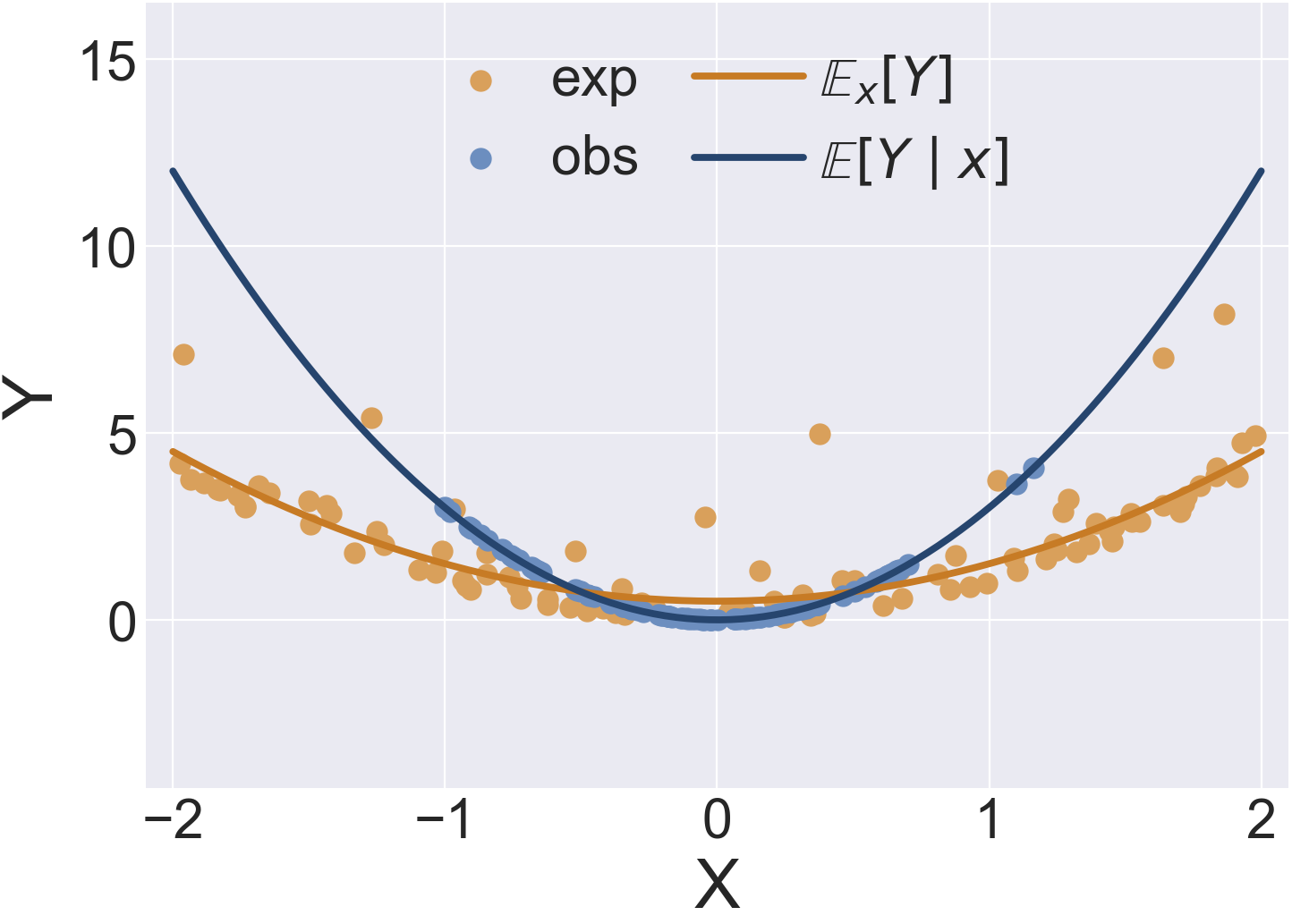}
		\caption{$\1M'$}
		\label{fig:unconfounded}
	\end{subfigure}\hfill\null
  \caption{Two confounded systems that induce the \emph{same} observational law $P(X, Y)$ (\cblue{blue}; regression $\mathbb{E}[Y \mid x]$ in \cdblue{dark blue}) yet \emph{different} interventional laws $\inv{Y}{x}$ (\corange{orange}), illustrating non-identifiability: (a) $\1M$ with $Y \gets -X^2 + U^2$ has a causal effect $\invE{Y}{x} = 1 - x^2$ (\cdorange{dark orange} curve); (b) $\1M'$ with $Y \gets X^2 + \tfrac{1}{2}U^2$ has a causal effect $\invE{Y}{x} = x^2 + \tfrac{1}{2}$ (\cdorange{dark orange} curve).}
  \label{fig:non-id}
\end{figure}
\begin{example}\label{exp:_1}
Consider a system $\1M$ with treatment $X \gets -U/2$ and outcome $Y \gets -X^2 + U^2$, where $U \sim \1N(0, 1)$ is unobserved. Fig.~\ref{fig:confounded} shows samples (blue) from the resulting observational distribution $P(X, Y)$ alongside samples (orange) from the interventional distribution $\inv{Y}{x}$ obtained by randomly assigning $X$ over $[-2, 2]$. The two deviate significantly: the unobserved confounder $U$ induces a spurious correlation between $X$ and $Y$, making some treatments appear more effective than they actually are. Treatment effects $\invE{Y}{x}$ are not identifiable here \citep{huang:val06,shpitser:pea06a}. To witness, consider the alternative system $\1M'$ with the same treatment $X \gets -U/2$ but outcome $Y \gets X^2 + \tfrac{1}{2} U^2$. Both models share the causal relationships ($X \rightarrow Y$ and $X \leftarrow U \rightarrow Y$) and induce the \emph{identical} observational law $P(X, Y)$ (substituting $U = -2X$ gives $Y = 3X^2$ in either model), yet they disagree on the effect: intervening $\doo(X \gets x)$ yields $\invE{Y}{x} = 1 - x^2$ under $\1M$ but $\invE{Y}{x} = x^2 + \tfrac{1}{2}$ under $\1M'$ (Fig.~\ref{fig:unconfounded}). Both models thus agree that $X$ affects $Y$ and are confounded in the same way; they disagree on \emph{which treatment is preferable}, with $x = 0$ optimal under $\1M$ but the worst choice under $\1M'$. The observational data alone cannot tell them apart. $\hfill \blacksquare$
\end{example}

When identifiability fails, a common alternative approach is \emph{partial identification}: bounding the effect over all confounding patterns compatible with the data \citep{manski1990nonparametric,balke:pea97}. For discrete treatments and outcomes, sharp bounds follow from programs over a discretized latent domain \citep{balke:pea97,zhang2022partial} or from generative-model approximations \citep{xia2021causal,xia2022neural}. However, this enumeration does not carry over to continuous treatments, where three lines of work exist.

\emph{Sensitivity models} bound how far hidden confounding may distort the treatment density: the continuous marginal sensitivity model (CMSM) of \citet{jesson2022scalable}, its sharp and generalized versions \citep{frauen2023sharp,baitairian2025sharp}, the $\delta$-MSM \citep{marmarelis2023partial}, and exposure-dependent and marginal-structural variants \citep{bonvini2022sensitivity,dalal2025partial}. Several bound the conditional average potential outcome (CAPO) pointwise, but all are stated relative to the observed treatment density and require it to be positive at the query. \emph{Optimization-based} bounds search a response-function class consistent with the data, usually given an instrument \citep{kilbertus2020class,gunsilius2020path,padh2023stochastic,balazadeh2022partial}. \emph{Smoothness} bounds restrict how fast counterfactual means vary \citep{tan2024consistency}. Following \citet{hausman2016individual}, \citet{kim2018identification} bound the effect of treatment on the treated assuming it is Lipschitz in the intervened treatment (smooth treatment response, STR) and in the received treatment (smooth treatment selection, STS), and \citet{khan2024offpolicy} use Lipschitz smoothness for off-policy evaluation without overlap. App.~\ref{appendix:related} provides more details on related work.

We work in the smoothness line, in the setting of \citet{kim2018identification} extended to covariates. Our targets are the CAPO $\tau(\*x, \*z) = \3E[Y_{\*x} \mid \*z]$ and the conditional average treatment effect (CATE) $\tau(\*x_1, \*z) - \tau(\*x_0, \*z)$ between two treatment levels. The bounds of \citet{kim2018identification} are sharp under STR together with a monotonicity assumption. Under STR and STS, however, sharp bounds remain an open challenge: we show that their combined bound is not sharp in general.

This paper aims to address this challenge. Our contributions are threefold. (1) \emph{Structural conditions for STR and STS.} We derive both conditions from Lipschitz continuity of the response functions and of the propensity scores given the latent confounder (Prop.~\ref{prop:structural}), with a selection cost that depends only on the observed treatment density (Lem.~\ref{lem:posterior_smooth}). (2) \emph{The bounds of \citet{kim2018identification} are loose.} Their combined STR+STS bound uses only two reference points on the identified diagonal (Prop.~\ref{prop:kim}) and is not sharp in general, even within their own model class (Prop.~\ref{prop:kimloose}). (3) \emph{Sharp bounds on the CAPO and the CATE.} Combining the STR and STS bounds through every diagonal reference point gives a valid bound (Thm.~\ref{thm:cate_bound}) that is the sharp identified set of the CAPO at every treatment, inside or outside the observational support (Thm.~\ref{thm:sharp}). Capping the difference of the upper and lower bounds gives the sharp identified set of the CATE, without any monotonicity assumption (Thm.~\ref{thm:contrast}). All proofs and experimental details are in the appendix.

\textbf{Notations.}
We denote variables by capital letters ($X$) and their values by lowercase ($x$); $|\*X|$ is the cardinality of a set $\*X$. $P(\*X)$ is the distribution of $\*X$, $P(\*x)$ abbreviates probability $P(\*X = \*x)$, and $P(\3X)$ stands for the probability $P(\*X \in \3X)$ that variables $\*X$ take a value in a set $\3X$. For reals $a, b$ we write $a \vee b \triangleq \max\{a, b\}$ and $a \wedge b \triangleq \min\{a, b\}$.

%% file: section2.tex
\section{Lipschitz Continuity for Partial Identification}\label{sec:_2}\label{sec:bounds}
Our analysis focuses on the structural causal models (SCMs) \citep{pearl:2k,bareinboim2022pearl} of a contextual bandit, with a treatment $\*X$, an outcome $Y$, covariates $\*Z$, and an unobserved confounder $\*U$ affecting all three. In a model $\1M$, structural functions set $\*Z \gets f_Z(\*U)$, $\*X \gets f_X(\*Z, \*U, \varepsilon_X)$, and $Y \gets f_Y(\*X, \*Z, \*U, \varepsilon_Y)$, where $\*U$, $\varepsilon_X$, and $\varepsilon_Y$ are mutually independent exogenous variables. The assignment noise $\varepsilon_X$ makes the treatment random given $(\*Z, \*U)$, so that it can have a conditional density $p(\*x \mid \*z, \*u)$; this is needed only for the structural condition \ref{ass:smoothness_X} below. The learner does not know $\1M$; it only observes samples of $(\*X, \*Z, Y)$, summarized by the \emph{observational distribution} $P(\*X, Y, \*Z) \triangleq P(\*X, Y, \*Z; \1M)$. An intervention $\doo(\*X \gets \*x)$ (for short, $\doo(\*x)$) replaces $f_X$ by the constant $\*x$, inducing the submodel $\1M_{\*x}$. The \emph{potential outcome} $Y_{\*x}(\*u)$ is the solution for $Y$ in $\1M_{\*x}$, and the \emph{interventional distribution} is the joint distribution $\inv{Y, \*Z}{\*x} \triangleq P\Parens{Y_{\*x}, \*Z_{\*x}; \1M}$ over potential outcomes.

Our goal is to infer the \emph{conditional average potential outcome} (CAPO) $\tau(\*x, \*z) \triangleq \invE{Y \mid \*z}{\*x}$ at every treatment $\*x$ and covariate value $\*z$ (\textbf{output}) from the observational distribution $P(\*X, \*Z, Y)$ (\textbf{input}). We also bound the \emph{conditional average treatment effect} (CATE) $\tau(\*x_1, \*z) - \tau(\*x_0, \*z) = \3E[Y_{\*x_1} - Y_{\*x_0} \mid \*z]$ between two treatment levels (Thm.~\ref{thm:contrast}). Under unobserved confounding, $\tau$ is not identifiable: many distinct SCMs reproduce the same $P(\*X, \*Z, Y)$ yet disagree on the target query. Consequently, we aim at partial identification: deriving informative bounds over the CAPO and CATO from the observational distribution.

Without restrictions on the unobserved mechanisms, the data leave $Y_{\*x}$ free at treatments that are rarely or never observed at $\*z$, and no estimator can improve on the trivial range bound there. An informative bound requires restricting how fast the unobserved mechanisms may vary \citep{kim2018identification,tan2024consistency}. We posit one for each unobserved mechanism: the outcome and the treatment assignment.

Each latent stratum induces a \emph{response function} $h_{\*u}(\*x, \*z) \triangleq \3E[f_Y(\*x, \*z, \*u, \varepsilon_Y)]$, the expected outcome of $\doo(\*X \gets \*x)$ at covariates $\*z$ within stratum $\*u$. Treatments and covariates take values in Borel sets $\3X \subseteq \3R^n$ and $\3Z \subseteq \3R^l$, which need not be bounded. With constants $L_Y, L_X > 0$, we assume the following.
\begin{enumerate}[label=(A\arabic*),start=1,leftmargin=0pt,itemindent=*,topsep=2pt,itemsep=1pt]
    \item \label{ass:bounded} (\emph{Compact outcome domain.}) The outcome takes values in a compact interval $\3Y = [\underline{y}, \overline{y}] \subset \3R$.
    \item \label{ass:smoothness_Y} (\emph{Lipschitz response functions.}) For all $\*u \in \3U$, $\*z \in \3Z$, and $\*x, \*x' \in \3X$,
    \begin{equation}\label{eq:holder_Y}
        |h_{\*u}(\*x, \*z) - h_{\*u}(\*x', \*z)| \le L_Y\, \|\*x - \*x'\|.
    \end{equation}
    \item \label{ass:smoothness_X} (\emph{Lipschitz propensity density.}) For all $\*u \in \3U$ and $\*z \in \3Z$, the treatment has a probability density $p(\*x \mid \*z, \*u)$ satisfying, for all $\*x, \*x' \in \3X$,
    \begin{equation}\label{eq:holder_X}
        \left \lvert p(\*x \mid \*z, \*u) - p(\*x' \mid \*z, \*u)\right \rvert \le L_X\, \|\*x - \*x'\|.
    \end{equation}
\end{enumerate}
We write $\operatorname{diam}(\3Y) \triangleq \overline{y} - \underline{y}$ for the outcome range. Every response function output $f_Y(\*x, \*z, \*u, \epsilon_Y)$ lies in $[\underline{y}, \overline{y}]$. The constants $L_Y$ and $L_X$ are \emph{sensitivity parameters} the analyst must prespecify.

\subsection{Smoothness of Treatment on the Treated}\label{sec:structural}
Assumptions~\ref{ass:smoothness_Y} and~\ref{ass:smoothness_X} restrict mechanisms that the data never reveal. To connect them to the data, we rewrite the target through a function whose values we partially identify. Specifically, let $g_{\*z}(\*x, \*x') \triangleq \3E[Y_{\*x} \mid \*x', \*z]$ denote the \emph{effect of treatment on the treated} (ETT) at treatment $\*x'$: the mean outcome under $\doo(\*x)$ of the units that observationally received $\*x'$ \citep{pearl:2k}. It is the function studied by \citet{kim2018identification}, here conditional on $\*z$. We call $\*x$ the \emph{intervened treatment} and $\*x'$ the \emph{received treatment}, and regard $g_{\*z}$ as a function on $\3X \times S_{\*z}$, where $S_{\*z} \triangleq \operatorname{supp} P(\*X \mid \*z)$ is the observational support at $\*z$. Conditioning on the received treatment $\*X = \*x'$ and applying the law of total expectation expands the target over the observational treatment distribution,
\begin{equation}\label{eq:cate_expand}
    \tau(\*x, \*z) = \int g_{\*z}(\*x, \*x')\, dP(\*x' \mid \*z).
\end{equation}
The \emph{nominal propensity score} $P(\*x' \mid \*z) = \int P(\*x' \mid \*z, \*u)\, dP(\*u \mid \*z)$ averages the latent-conditional assignments over the strata, and we write $p(\*x' \mid \*z)$ for its density. As the nominal propensity is identified, bounding $\tau$ reduces to bounding $g_{\*z}$. On the diagonal $\*x' = \*x$ the ETT is point-identified by consistency \citep[Ch.~7]{pearl:2k}: $g_{\*z}(\*x, \*x) = \mu(\*x, \*z)$ for $\*x \in S_{\*z}$, where $\mu(\*x, \*z) \triangleq \3E[Y \mid \*x, \*z]$ is the observed regression. Off the diagonal, where $\*x \neq \*x'$, it is unidentified.

Assumptions~\ref{ass:smoothness_Y} and~\ref{ass:smoothness_X} restrict the unobserved mechanisms, but the bounds use them only through two consequences for the ETT function, which, following \citet{kim2018identification}, limit how fast it may vary in each coordinate. Note that $Y_{\*x} \perp \*X \mid \*U, \*Z$ in every SCM of $\2M$, and the received treatment enters the ETT only through the latent posterior $P(\*U \mid \*x', \*z)$. Expanding the total expectation over $\*U$ gives
\begin{equation}\label{eq:ett}
    g_{\*z}(\*x, \*x') = \int h_{\*u}(\*x, \*z)\, dP(\*u \mid \*x', \*z).
\end{equation}
This splits each ETT value into a response function, which varies with the intervened treatment and is controlled by \ref{ass:smoothness_Y}, and a latent posterior, which varies with the received treatment and is controlled by \ref{ass:smoothness_X}. The key step bounds how fast the latent posterior moves with the received treatment, using only the density $p(\*x \mid \*z)$ of the nominal propensity.
\begin{restatable}{lemma}{lemposteriorsmooth}\label{lem:posterior_smooth}
Under \ref{ass:smoothness_X}, for every $\*z \in \3Z$ and every pair $\*x_0', \*x_1' \in \3X$ with $p(\*x_0' \mid \*z) > 0$ and $p(\*x_1' \mid \*z) > 0$,
\begin{equation}\label{eq:posterior_smooth}
    \mathrm{TV}_{\*u}\bigl(\*x_0', \*x_1'; \*z\bigr) \le \frac{L_X\,\|\*x_0' - \*x_1'\|}{\max\{p(\*x_0' \mid \*z),\, p(\*x_1' \mid \*z)\}},
\end{equation}
where $\mathrm{TV}_{\*u}(\*x_0', \*x_1'; \*z) \triangleq \mathrm{TV}\bigl(P(\*U \mid \*x_0', \*z),\, P(\*U \mid \*x_1', \*z)\bigr)$ is the total variation between latent posteriors.
\end{restatable}
Two received treatments share the response functions of Eq.~\ref{eq:ett} and differ only through the latent posterior, and the response functions take values in $\3Y$. So, for every $\*x \in \3X$ and every such pair,
\begin{equation}\label{eq:transport_xp}
\begin{aligned}
    &\bigl|g_{\*z}(\*x, \*x_0') - g_{\*z}(\*x, \*x_1')\bigr| \le \operatorname{diam}(\3Y)\,\mathrm{TV}_{\*u}(\*x_0', \*x_1'; \*z) \\
    &\quad \le \kappa_{\*z}(\*x_0', \*x_1') \triangleq \frac{L_X\operatorname{diam}(\3Y)\,\|\*x_0' - \*x_1'\|}{\max\{p(\*x_0' \mid \*z),\, p(\*x_1' \mid \*z)\}}.
\end{aligned}
\end{equation}
The coefficient $\kappa_{\*z}$ need not satisfy the triangle inequality, so we pass to its chain closure $\kappa^{\ast}_{\*z}(\*a, \*b) \triangleq \inf \sum_{i} \kappa_{\*z}(\*a_{i}, \*a_{i+1})$ over finite chains $\*a = \*a_0, \dots, \*a_k = \*b$ of points with positive density. This coefficient is \emph{density-adaptive}: it is small where many units are observed and grows where the observed density thins out, and it evaluates the density only at support points. A function satisfying $|f(\*a) - f(\*b)| \le \kappa_{\*z}(\*a, \*b)$ for all pairs also satisfies it with $\kappa^{\ast}_{\*z}$, so the constraint translates.
\begin{restatable}{proposition}{propstructural}\label{prop:structural}
Under \ref{ass:bounded}, \ref{ass:smoothness_Y}, and \ref{ass:smoothness_X}, for every $\*z \in \3Z$ the ETT function satisfies the following two conditions with the costs $\rho_{\mathrm{T}} \triangleq L_Y\|\cdot\|$ on $\3X$ and $\rho_{\mathrm{S}} \triangleq \kappa^{\ast}_{\*z}$ on $S^{+}_{\*z} \triangleq \{\*x' \in S_{\*z} : p(\*x' \mid \*z) > 0\}$.
\begin{enumerate}[label=(STR),leftmargin=0pt,itemindent=*,topsep=2pt,itemsep=1pt]
    \item \label{ass:str} (\emph{Smooth treatment response.}) For all $\*x, \tilde{\*x} \in \3X$ and $\*x' \in S^{+}_{\*z}$,
    \begin{equation}\label{eq:str}
        |g_{\*z}(\*x, \*x') - g_{\*z}(\tilde{\*x}, \*x')| \le \rho_{\mathrm{T}}(\*x, \tilde{\*x}).
    \end{equation}
\end{enumerate}
\begin{enumerate}[label=(STS),leftmargin=0pt,itemindent=*,topsep=2pt,itemsep=1pt]
    \item \label{ass:sts} (\emph{Smooth treatment selection.}) For all $\*x \in \3X$ and $\*x_0', \*x_1' \in S^{+}_{\*z}$,
    \begin{equation}\label{eq:sts}
        |g_{\*z}(\*x, \*x_0') - g_{\*z}(\*x, \*x_1')| \le \rho_{\mathrm{S}}(\*x_0', \*x_1').
    \end{equation}
\end{enumerate}
\end{restatable}
\ref{ass:str} bounds how fast the counterfactual mean of a fixed group of units changes with the intervened treatment. \ref{ass:sts} bounds how fast it changes across the groups selected by different observed treatments, which is where confounding enters: without confounding $g_{\*z}$ is constant in $\*x'$, and \ref{ass:sts} holds with $\rho_{\mathrm{S}} = 0$. These are the smooth treatment response and smooth treatment selection conditions of \citet{kim2018identification}, stated conditionally on $\*z$.

For each covariate value $\*z$, let $\rho_{\mathrm{T}}$ be a pseudometric on $\3X$ and $\rho_{\mathrm{S}}$ a pseudometric on $S_{\*z}$, both continuous in the Euclidean topology (we suppress their dependence on $\*z$), and say that the ETT function satisfies \ref{ass:str} and \ref{ass:sts} with these costs when the two inequalities of Prop.~\ref{prop:structural} hold with them for all $\*x, \tilde{\*x} \in \3X$ and $\*x', \*x_0', \*x_1' \in S_{\*z}$. We call $\rho_{\mathrm{T}}$ the \emph{response cost} and $\rho_{\mathrm{S}}$ the \emph{selection cost}. The default choices $\rho_{\mathrm{T}} = L_Y\|\cdot\|$ and $\rho_{\mathrm{S}} = C\|\cdot\|$, with constants $L_Y, C > 0$, are the constants $b$ and $a$ of \citet{kim2018identification}; a general $\rho_{\mathrm{T}}$ extends their modulus of continuity, and a general $\rho_{\mathrm{S}}$ admits the density-adaptive $\kappa^{\ast}_{\*z}$. Like $L_Y$ and $L_X$, the costs are sensitivity parameters.

The rest of the paper proceeds as follows. Sec.~\ref{sec:cate_bounds} derives the STR and STS bounds and combines them through every diagonal reference (Thm.~\ref{thm:cate_bound}). Sec.~\ref{sec:sharpness} shows that the combined bound of \citet{kim2018identification} is a special case of ours and is not sharp in general, whereas ours is the sharp identified set of the CAPO and of the CATE (Sec.~\ref{sec:sharpset}). The bounding results of Sec.~\ref{sec:cate_bounds} are derived under \ref{ass:bounded}--\ref{ass:smoothness_X} with the costs of Prop.~\ref{prop:structural}, but extend to \ref{ass:bounded}, \ref{ass:str}, and \ref{ass:sts} with any costs $\rho_{\mathrm{T}}, \rho_{\mathrm{S}}$ as defined above (App.~\ref{appendix:a}).

%% file: section3.tex
\section{Lipschitz Bounds on Treatment Responses}\label{sec:holder}\label{sec:cate_bounds}
\begin{figure*}[t]
\centering
\includegraphics[width=\textwidth]{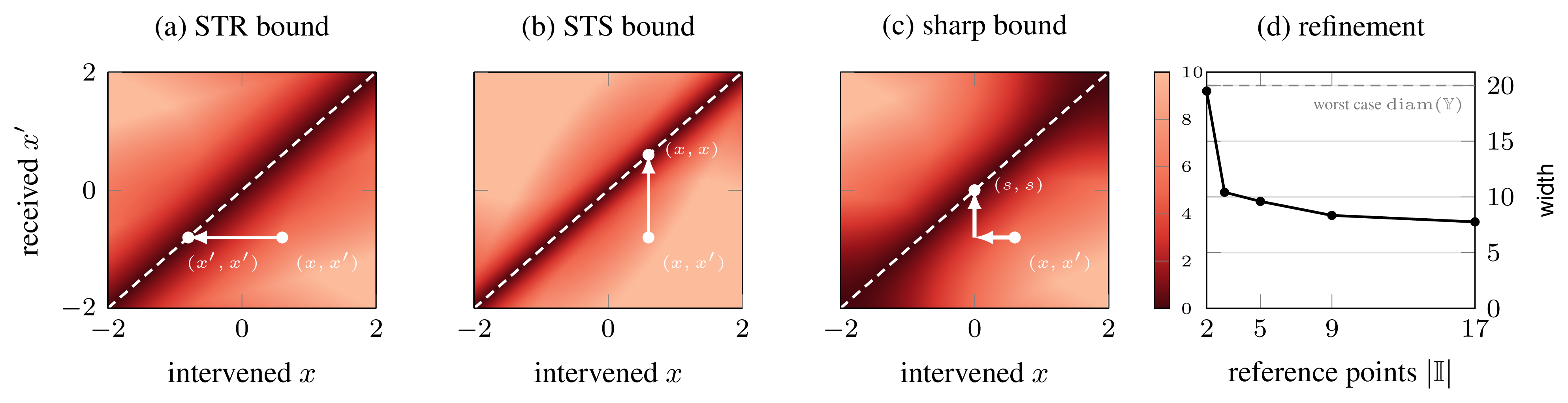}
\caption{Half-width of the bound on the ETT over the $(x,x')$ plane for Exs.~\ref{exp:dox} and~\ref{exp:bounds}. The identified diagonal $x=x'$ is dashed, and arrows point from each query to the diagonal point bounding it. (a) STR bound (\ref{ass:str}); (b) STS bound (\ref{ass:sts}); (c) sharp bound over all diagonal reference points (Thm.~\ref{thm:sharp}); (d) enlarging a finite reference set $\3I$ shrinks the averaged width toward the sharp bound. Panels (a)--(c) saturate at $\tfrac{1}{2}\operatorname{diam}(\3Y)=10$.}
\label{fig:bound_heatmap}
\end{figure*}

By Eq.~\ref{eq:cate_expand}, bounding $\tau(\*x, \*z)$ reduces to bounding the ETT $g_{\*z}(\*x, \*x')$, of which only the diagonal is identifiable. We first derive the STR and STS bounds of \citet{kim2018identification}, conditional on $\*z$, under structural conditions; each moves one argument of $g_{\*z}$ onto the diagonal. Sec.~\ref{sec:combined} combines them to derive the sharp bound, over all possible pathways to transport an unidentified argument to an identified one on the diagonal.

By consistency, \ref{ass:str} places each $g_{\*z}(\*x, \*x')$ within $\rho_{\mathrm{T}}(\*x, \*x')$ of $\mu(\*x', \*z) = g_{\*z}(\*x', \*x')$ (Fig.~\ref{fig:bound_heatmap}(a)).

\begin{restatable}{lemma}{lemsidebound}\label{lem:side_bound}
Under \ref{ass:bounded}--\ref{ass:smoothness_X}, for every $(\*x, \*z) \in \3X \times \3Z$, the following inequalities hold:
\begin{equation}\label{eq:side_anchor}
    \underline{w}(\*x, \*z) \;\le\; \3E[Y_{\*x} \mid \*z] \;\le\; \overline{w}(\*x, \*z),
\end{equation}
where, for a draw $\*X' \sim P(\*X \mid \*z)$ of the received treatment, the bounds are
\begin{equation}\label{eq:side_halfwidth}
\begin{aligned}
    \underline{w}(\*x, \*z) &\triangleq \3E\bigl[\,\bigl(\mu(\*X', \*z) - \rho_{\mathrm{T}}(\*x, \*X')\bigr) \vee \underline{y} \,\bigm|\, \*z\bigr], \\
    \overline{w}(\*x, \*z) &\triangleq \3E\bigl[\,\bigl(\mu(\*X', \*z) + \rho_{\mathrm{T}}(\*x, \*X')\bigr) \wedge \overline{y} \,\bigm|\, \*z\bigr].
\end{aligned}
\end{equation}
\end{restatable}
Lem.~\ref{lem:side_bound} is \citet[Prop.~3.1]{kim2018identification}, conditional on $\*z$ and intersected with $\3Y$. It uses no information at the query, so it holds at every $\*x \in \3X$, inside or outside $S_{\*z}$; without confounding it sharpens to the bound of \citet{khan2024offpolicy} (Cor.~\ref{cor:khan}).

\ref{ass:sts} instead moves the received treatment, placing $g_{\*z}(\*x, \*x')$ within $\rho_{\mathrm{S}}(\*x', \*x)$ of $\mu(\*x, \*z)$ (Fig.~\ref{fig:bound_heatmap}(b)), which requires the query to be observed.
\begin{restatable}{lemma}{lemcondbound}\label{lem:cond_bound}
Under \ref{ass:bounded}--\ref{ass:smoothness_X}, for every $\*z \in \3Z$ and $\*x \in S_{\*z}$, the following inequalities hold:
\begin{equation}\label{eq:cond_side_anchor}
    \underline{w}_X(\*x, \*z) \;\le\; \3E[Y_{\*x} \mid \*z] \;\le\; \overline{w}_X(\*x, \*z),
\end{equation}
where, for a draw $\*X' \sim P(\*X \mid \*z)$ of the received treatment, the bounds are
\begin{equation}\label{eq:cond_halfwidth}
\begin{aligned}
    \underline{w}_X(\*x, \*z) &\triangleq \3E\bigl[\,\bigl(\mu(\*x, \*z) - \rho_{\mathrm{S}}(\*X', \*x)\bigr) \vee \underline{y} \,\bigm|\, \*z\bigr], \\
    \overline{w}_X(\*x, \*z) &\triangleq \3E\bigl[\,\bigl(\mu(\*x, \*z) + \rho_{\mathrm{S}}(\*X', \*x)\bigr) \wedge \overline{y} \,\bigm|\, \*z\bigr].
\end{aligned}
\end{equation}
\end{restatable}
Lem.~\ref{lem:cond_bound} is \citet[Prop.~D.1]{kim2018identification}, conditional on a context $\*z$. With  costs measured by Euclidean distance between treatments, it beats the STR bound wherever $C < L_Y$. With the density-adaptive cost of Prop.~\ref{prop:structural} it needs no overlap floor: its averaged band is at most $\3E[\kappa_{\*z}(\*x, \*X') \mid \*z] \le L_X\operatorname{diam}(\3Y)\int_{S_{\*z}}\|\*x - \*x'\|\, d\*x'$, since $\max\{p(\*x \mid \*z), p(\*x' \mid \*z)\} \ge p(\*x' \mid \*z)$, which is finite for bounded $S_{\*z}$ at any $\*x \in S^{+}_{\*z}$, however small $p(\*x \mid \*z)$ is.

\begin{example}[STR and STS bounds]\label{exp:dox}
For the SCM of Ex.~\ref{exp:_1} ($X \gets -U/2$, $Y \gets -X^2 + U^2$), \ref{ass:bounded} requires bounded outcomes, so we truncate $U$ to $[-4, 4]$, which puts the treatments in $\3X = [-2, 2]$ and the outcomes in $\3Y = [-4, 16]$. The response functions are $h_{u}(x) = -x^2 + u^2$, which sweep exactly $\3Y$ over that region, so $|h_{u}(x) - h_{u}(x')| = |x + x'|\,|x - x'| \le 4|x - x'|$ on $[-2, 2]$: the outcome mechanism is $4$-Lipschitz, and \ref{ass:str} holds with $\rho_{\mathrm{T}} = 4|\cdot|$ (Prop.~\ref{prop:structural}). The observed treatments fill $\3X$, and along that manifold $\mu(x) = 3x^2 \in [0, 12] \subset \3Y$, so the diagonal $\3E[Y_x \mid x] = \mu(x)$ is identified by consistency. The STR bound (Lem.~\ref{lem:side_bound}) carries these identified diagonal values to a query $x$, averaging over $P(X \mid \*z)$ the interval $[3X'^2 - 4|x - X'|,\ 3X'^2 + 4|x - X'|] \cap [-4, 16]$. Assumption~\ref{ass:smoothness_X} cannot hold in this system, since $X = -U/2$ is a deterministic function of $U$, so we posit \ref{ass:sts} directly, with $\rho_{\mathrm{S}} = C|\cdot|$ and $C = 8$. The STS bound (Lem.~\ref{lem:cond_bound}) instead centers the ETT at the diagonal value at the \emph{query} treatment, $\mu(x) = 3x^2$, averaging the interval $[3x^2 - 8|x - X'|,\ 3x^2 + 8|x - X'|] \cap [-4, 16]$. Since $C = 8 > L_Y = 4$, here the STR bound is tighter. $\hfill \blacksquare$
\end{example}

\subsection{Combined Bounds over CAPOs}\label{sec:combined}
Each of the STR and STS bounds moves a single argument of the unidentified ETT $g_{\*z}(\*x, \*x')$ onto the identified diagonal: the intervened treatment $\*x \mapsto \*x'$ in the first case, the received treatment $\*x' \mapsto \*x$ in the second. Applying both, through the intermediate point $(\*s, \*x')$, bounds the ETT by the identified value $\mu(\*s, \*z)$ at \emph{any} diagonal reference point $\*s \in S_{\*z}$: a step $\*x \to \*s$ in the intervened treatment at cost $\rho_{\mathrm{T}}(\*x, \*s)$, then a step $\*x' \to \*s$ in the received treatment at cost $\rho_{\mathrm{S}}(\*x', \*s)$. Formally, the two smoothness conditions are one Lipschitz condition for the product cost $d_{\*z}((\*x, \*x'), (\*a, \*b)) \triangleq \rho_{\mathrm{T}}(\*x, \*a) + \rho_{\mathrm{S}}(\*x', \*b)$ on $\3X \times S_{\*z}$.

\begin{restatable}{lemma}{lemrectanchor}\label{lem:rect_anchor}
Under \ref{ass:bounded}--\ref{ass:smoothness_X}, for every $\*z \in \3Z$ the ETT function is Lipschitz for the product cost: $|g_{\*z}(\*p) - g_{\*z}(\*q)| \le d_{\*z}(\*p, \*q)$ for all $\*p, \*q \in \3X \times S_{\*z}$. In particular, for every $\*x \in \3X$ and $\*x', \*s \in S_{\*z}$,
\begin{align}
    &\bigl|\3E[Y_{\*x} \mid \*x', \*z] - \mu(\*s, \*z)\bigr| \le \rho(\*s, \*x') \label{eq:rect_anchor} \\
    &\text{with } \rho(\*s, \*x') \triangleq \rho_{\mathrm{T}}(\*x, \*s) + \rho_{\mathrm{S}}(\*x', \*s). \label{eq:anchor_radius}
\end{align}
\end{restatable}
Lem.~\ref{lem:rect_anchor} contains both STR and STS bounds: $\*s = \*x'$ removes the STS step and recovers the pointwise STR bound of Lem.~\ref{lem:side_bound}, while $\*s = \*x$ removes the STR step and recovers the pointwise STS bound of Lem.~\ref{lem:cond_bound}. Taking the better of these two at each $\*x'$ is the combined STR+STS bound of \citet[Prop.~D.3]{kim2018identification}.

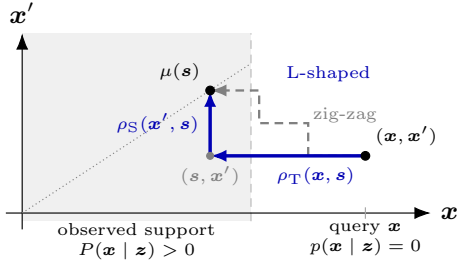
\begin{figure}[t]
\centering
\resizebox{0.75\linewidth}{!}{%
\begin{tikzpicture}[every node/.style={font=\footnotesize}]
  \fill[gray!12] (-0.8,-0.1) rectangle (2.0,2.2);
  \draw[gray!55,densely dashed] (2.0,-0.1) -- (2.0,2.2);
  \draw[->] (-1.0,0) -- (4.2,0) node[right]{$\*x$};
  \draw[->] (-0.8,-0.2) -- (-0.8,2.2) node[above]{$\*x'$};
  \draw[gray,densely dotted] (-0.8,0) -- (2.0,1.83);
  \coordinate (P) at (3.4,0.7); \coordinate (A) at (1.5,1.5); \coordinate (Lc) at (1.5,0.7);
  \draw[->,gray,thick,densely dashed] (P)--(2.7,0.7)--(2.7,1.1)--(2.1,1.1)--(2.1,1.5)--(A);
  \draw[->,blue!70!black,very thick] (P)--(Lc) node[pos=0.32,below,font=\tiny]{$\rho_{\mathrm{T}}(\*x,\*s)$};
  \draw[->,blue!70!black,very thick] (Lc)--(A) node[midway,left,font=\tiny]{$\rho_{\mathrm{S}}(\*x',\*s)$};
  \filldraw (P) circle(1.5pt) node[above right,font=\tiny]{$(\*x,\*x')$};
  \filldraw (A) circle(1.5pt) node[above left,font=\tiny]{$\mu(\*s)$};
  \filldraw[gray] (Lc) circle(1.2pt) node[below,gray,font=\tiny]{$(\*s,\*x')$};
  \node[blue!70!black,font=\tiny] at (2.95,1.7) {L-shaped};
  \node[gray,font=\tiny] at (3.15,1.2) {zig-zag};
  \draw[gray!70] (3.4,-0.07) -- (3.4,0.07);
  \node[font=\tiny,align=center,anchor=north] at (0.6,0) {observed support\\$P(\*x\mid\*z) > 0$};
  \node[font=\tiny,align=center,anchor=north] at (3.4,0) {query $\*x$\\$p(\*x\mid\*z)=0$};
\end{tikzpicture}}
\caption{The L-shaped path of Lem.~\ref{lem:rect_anchor} (blue) from $(\*x,\*x')$ through $(\*s,\*x')$ to the diagonal point $(\*s,\*s)$, at cost $\rho(\*s,\*x')$ (Eq.~\ref{eq:anchor_radius}). A zig-zag path (grey) costs at least as much; query $\*x$ may lie outside the support.}
\label{fig:noint_L}
\end{figure}

Lem.~\ref{lem:rect_anchor} reaches the reference point along an L-shaped path, one step in each argument (Fig.~\ref{fig:noint_L}, blue). By the triangle inequality for each cost, no path that alternates the two arguments, and no path through several diagonal points, is cheaper (App.~\ref{appendix:a}), and Thm.~\ref{thm:sharp} shows more: no argument at all can improve the resulting bound, because it is realizable by a candidate SCM compatibel with the observational distribution.

Taking every diagonal point as a reference point gives, for each covariate value $\*z$, the \emph{upper and lower bounds}
\begin{equation}\label{eq:envelopes}
\begin{aligned}
    U_{\*z}(\*x, \*x') &\triangleq \Bigl(\inf_{\*s \in S_{\*z}} \bigl[\mu(\*s, \*z) + \rho(\*s, \*x')\bigr]\Bigr) \wedge \overline{y}, \\
    L_{\*z}(\*x, \*x') &\triangleq \Bigl(\sup_{\*s \in S_{\*z}} \bigl[\mu(\*s, \*z) - \rho(\*s, \*x')\bigr]\Bigr) \vee \underline{y},
\end{aligned}
\end{equation}
with $\rho$ as in Eq.~\ref{eq:anchor_radius}. Since $\rho(\*s, \*x')$ is the product cost $d_{\*z}\bigl((\*x, \*x'), (\*s, \*s)\bigr)$ from the query to the diagonal point $(\*s, \*s)$, $U_{\*z}$ and $L_{\*z}$ are the McShane--Whitney extensions of the diagonal values $\mu(\*s, \*z)$ to $\3X \times S_{\*z}$ \citep{mcshane1934extension,whitney1934analytic}, truncated to $\3Y$; Sec.~\ref{sec:sharpset} shows that these bounds are extremal.
\begin{restatable}{theorem}{thmcatebound}\label{thm:cate_bound}
Under \ref{ass:bounded}--\ref{ass:smoothness_X}, for every $\*z \in \3Z$ the ETT function satisfies $L_{\*z} \le g_{\*z} \le U_{\*z}$ on $\3X \times S_{\*z}$. Consequently, for every $\*x \in \3X$, including $\*x \notin S_{\*z}$, and a draw $\*X' \sim P(\*X \mid \*z)$ of the received treatment,
\begin{equation}\label{eq:cate_valid}
    \3E\bigl[L_{\*z}(\*x, \*X') \mid \*z\bigr] \le \tau(\*x, \*z) \le \3E\bigl[U_{\*z}(\*x, \*X') \mid \*z\bigr].
\end{equation}
\end{restatable}
For a query $\*x \notin S_{\*z}$, the regression $\mu(\*x, \*z)$ is undefined, so the STS bound of Lem.~\ref{lem:cond_bound} is unavailable, and the combined bound of \citet[Prop.~D.3]{kim2018identification} reduces to STR alone. The bound of Eq.~\ref{eq:envelopes} still uses both steps, through reference points $\*s \in S_{\*z}$ near $\*x$ that trade $\rho_{\mathrm{T}}(\*x, \*s)$ against $\rho_{\mathrm{S}}(\*x', \*s)$. With the density-adaptive cost of Prop.~\ref{prop:structural}, the selection step runs only between support points, so it stays informative even when little mass lies near the query, i.e., the bound needs no positivity assumption at the query.

Thm.~\ref{thm:cate_bound} is Lem.~\ref{lem:rect_anchor} applied at every reference point, together with the range constraint. It holds for every SCM satisfying \ref{ass:bounded}--\ref{ass:smoothness_X}, and more generally for every SCM whose ETT function satisfies \ref{ass:str} and \ref{ass:sts} with arbitrary costs. Over that class, these bounds cannot be improved (Sec.~\ref{sec:sharpset}).

In practice, we compute the bounds over a finite set of reference points in the support, where $\mu$ is identified; restricting the references keeps the bounds valid.
\begin{restatable}{corollary}{corfinite}\label{cor:finite}
Under \ref{ass:bounded}--\ref{ass:smoothness_X}, let $\3I \subseteq S_{\*z}$ be a nonempty reference set, and let $U_{\*z, \3I}$ and $L_{\*z, \3I}$ be the bounds of Eq.~\ref{eq:envelopes} with the infimum and supremum restricted to $\3I$. Then $L_{\*z, \3I} \le L_{\*z} \le U_{\*z} \le U_{\*z, \3I}$, so that
\begin{equation}\label{eq:finite_bound}
    \3E\bigl[L_{\*z, \3I}(\*x, \*X') \mid \*z\bigr] \le \tau(\*x, \*z) \le \3E\bigl[U_{\*z, \3I}(\*x, \*X') \mid \*z\bigr]
\end{equation}
for every $\*x \in \3X$. The bounds tighten as $\3I$ grows. For finite $\3I$, they cost $O(n|\3I|)$ per query when $P(\*X \mid \*z)$ is an empirical law on $n$ points.
\end{restatable}
Since the STR and STS bounds use the references $\*s = \*x'$ and $\*s = \*x$, the combined bound is never looser than either; Cor.~\ref{cor:strict} in App.~\ref{appendix:a} states when it is tighter.

\begin{example}[Combined bound]\label{exp:bounds}
Continuing Ex.~\ref{exp:dox} ($\rho_{\mathrm{T}} = 4|\cdot|$, $\rho_{\mathrm{S}} = 8|\cdot|$, $\mu(s) = 3s^2$ on $S = [-2, 2]$). The bound of Thm.~\ref{thm:cate_bound} takes the best diagonal reference for each pair $(x, x')$ (Fig.~\ref{fig:bound_heatmap}(c)), and an intermediate point can beat both the STR and the STS bound. At $(x, x') = (-1, 1)$, the STR bound gives the upper value $\mu(1) + 4 \cdot 2 = 11$ and the STS bound $\mu(-1) + 8 \cdot 2 = 19$, whereas the reference $s = 2/3$ gives $\tfrac{4}{3} + \tfrac{20}{3} + \tfrac{8}{3} = \tfrac{32}{3} \approx 10.67$. Averaging over $P(X \mid \*z)$ gives the sharp bound over true response $\tau(x) = \3E[U^2] - x^2 \approx 1 - x^2$; finite references (Cor.~\ref{cor:finite}) give widths of Fig.~\ref{fig:bound_heatmap}(d), which shrink as $\3I$ grows. $\hfill \blacksquare$
\end{example}

%% file: section4.tex
\section{Sharpness Analysis}\label{sec:sharpness}
Thm.~\ref{thm:cate_bound} shows that the bound is valid. This section shows that it cannot be improved. Sharpness is relative to a model class, and we establish it for the class defined by \ref{ass:str} and \ref{ass:sts} with the general costs of Sec.~\ref{sec:structural}. Let $\2M(\rho_{\mathrm{T}}, \rho_{\mathrm{S}})$ be the class of SCMs satisfying \ref{ass:bounded} whose ETT functions satisfy \ref{ass:str} and \ref{ass:sts}; by Prop.~\ref{prop:structural}, it contains every SCM satisfying \ref{ass:bounded}--\ref{ass:smoothness_X} with the costs of that proposition. Given the observational law, the \emph{identified set} of $\tau(\*x, \*z)$ is
\begin{equation}\label{eq:identified_set}
\begin{split}
    &\Theta(\*x, \*z) \triangleq \bigl\{\, \3E_{\1M}[Y_{\*x} \mid \*z] :\ \1M \in \2M(\rho_{\mathrm{T}}, \rho_{\mathrm{S}}), \\
    &\qquad P_{\1M}(\*X, \*Z, Y) = P(\*X, \*Z, Y) \,\bigr\}.
\end{split}
\end{equation}
We first compare our bound with the combined bound of \citet{kim2018identification}, which uses only two of its reference points and can be tightened. We then show that our bound cannot: it gives $\Theta(\*x, \*z)$ exactly, and it yields the sharp identified set of the CATE (Sec.~\ref{sec:sharpset}).

\citet{kim2018identification} combine the two conditions by taking, at each received treatment, the better of the STR and STS bounds of Sec.~\ref{sec:cate_bounds}. This restricts the bounds in Eq.~\ref{eq:envelopes} to two reference points.
\begin{restatable}{proposition}{propkim}\label{prop:kim}
Let $\3X \subset \3R$, without covariates, and $\rho_{\mathrm{T}} = b|\cdot|$, $\rho_{\mathrm{S}} = a|\cdot|$. For any query $x \in S$, the combined STR+STS upper bound of \citet[Prop.~D.3]{kim2018identification} gives the following
\begin{equation*}
    u_2(x) \triangleq \3E\bigl[\min\bigl\{\mu(x) + a|X' - x|,\ \mu(X') + b|X' - x|\bigr\}\bigr],
\end{equation*}
which is the bound of Cor.~\ref{cor:finite} with the two-point reference set $\{x, x'\}$ chosen separately at each received treatment $x'$ (without the intersection with $\3Y$).
\end{restatable}
Dropping the other reference points costs sharpness.
\begin{restatable}{proposition}{propkimloose}\label{prop:kimloose}
The bound $u_2$ of Prop.~\ref{prop:kim} is not sharp in general: for some costs, observational distribution, and query $x \in S$, $u_2(x)$ strictly exceeds $\3E[Y_x]$ in every SCM satisfying the STR and STS conditions with these costs and reproducing the observational distribution.
\end{restatable}
For example, let $a = b = 1$, let the treatment be uniform on $S = \{0, 1, 2\}$, let $Y = \mu(X)$ with $\mu = (0, -1, 0)$, and take the query $x = 2$. Table~\ref{tab:kim} lists the terms $\mu(s) + |2 - s| + |x' - s|$. The two-point bound uses only the columns $s \in \{2, x'\}$, which give the values $(2, 0, 0)$ across $x' = 0, 1, 2$ and hence $u_2(x = 2) = 2/3$. The reference point $s = 1$ lowers the value at $x' = 0$ from $2$ to $1$: every such model has $\3E[Y_{x=2}] \le 1/3$, and App.~\ref{appendix:a} constructs one that attains this upper bound. The sharpness argument for Prop.~D.3 in \citet{kim2018identification} constructs a data-generating process whose regression $\mu$ is constant. With a constant regression and these costs, $\inf_s [b|x - s| + a|x' - s|] = \min\{a, b\}|x - x'|$ is attained at $s \in \{x, x'\}$, so the two-point bound coincides with our bound; for a nonconstant regression, as here, it does not and can be generally improved.
\begin{table}[t]
\centering
\small
\begin{tabular}{c|ccc|cc}
$x'$ & $s = 0$ & $s = 1$ & $s = 2$ & two-point & sharp \\ \hline
$0$ & $2$ & $1$ & $2$ & $2$ & $1$ \\
$1$ & $3$ & $0$ & $1$ & $0$ & $0$ \\
$2$ & $4$ & $1$ & $0$ & $0$ & $0$
\end{tabular}
\caption{The counterexample of Prop.~\ref{prop:kimloose} at the query $x = 2$: terms $\mu(s) + |2 - s| + |x' - s|$, the two-point bound of \citet{kim2018identification} (columns $s \in \{2, x'\}$), and the sharp bound. Averaging over $P$ gives $2/3$ and $1/3$.}
\label{tab:kim}
\end{table}

\subsection{Sharp Identified Set}\label{sec:sharpset}
The diagonal of any function satisfying \ref{ass:str} and \ref{ass:sts} is Lipschitz for $\rho_{\mathrm{T}} + \rho_{\mathrm{S}}$: for all $\*s, \*s' \in S_{\*z}$,
\begin{equation}\label{eq:consistency}
    |\mu(\*s, \*z) - \mu(\*s', \*z)| \le \rho_{\mathrm{T}}(\*s, \*s') + \rho_{\mathrm{S}}(\*s, \*s'),
\end{equation}
and in particular it is continuous. This is a \emph{testable implication} of the model, and it is the only restriction the model places on the observational law (Cor.~\ref{cor:test} in App.~\ref{appendix:a}). When it holds, Thm.~\ref{thm:sharp} characterizes $\Theta(\*x, \*z)$ exactly at every query, including treatments $\*x$ outside the observational support at $\*z$.
\begin{restatable}{theorem}{thmsharp}\label{thm:sharp}
Under \ref{ass:bounded}, \ref{ass:str}, and \ref{ass:sts}, fix $\*z \in \3Z$ and suppose $\mu(\cdot, \*z)$ satisfies Eq.~\ref{eq:consistency}. Then:
\begin{enumerate}[label=(\alph*),leftmargin=*,topsep=2pt,itemsep=1pt]
    \item (\emph{Attainment}) $U_{\*z}$ and $L_{\*z}$ take values in $\3Y$, satisfy the inequalities of \ref{ass:str} and \ref{ass:sts}, and equal $\mu(\*s, \*z)$ at every diagonal point $(\*s, \*s)$;
    \item (\emph{Sharp set}) for every $\*x \in \3X$, including $\*x \notin S_{\*z}$,
    \begin{equation}\label{eq:cate_two_sided}
         \Theta(\*x, \*z) = \Bigl[\,\3E\bigl[L_{\*z}(\*x, \*X') \mid \*z\bigr], \3E\bigl[U_{\*z}(\*x, \*X') \mid \*z\bigr]\,\Bigr].
    \end{equation}
\end{enumerate}
If Eq.~\ref{eq:consistency} fails on a set of covariate values of positive probability, no SCM in $\2M(\rho_{\mathrm{T}}, \rho_{\mathrm{S}})$ reproduces $P(\*X, \*Z, Y)$, and the model is refuted.
\end{restatable}
Part (a) holds because an infimum of $d_{\*z}$-Lipschitz functions is $d_{\*z}$-Lipschitz, while Eq.~\ref{eq:consistency} pins the bounds to the data on the diagonal. So the McShane--Whitney extensions $U_{\*z}$ and $L_{\*z}$ of Sec.~\ref{sec:combined} are the largest and the smallest functions with values in $\3Y$ that satisfy \ref{ass:str} and \ref{ass:sts} and agree with $\mu$ on the diagonal. Part (b) constructs, for every admissible ETT function, an SCM that realizes it while reproducing the observational distribution; mixing the two bounds then sweeps the whole interval (App.~\ref{appendix:a}).

The CATE between two treatment levels, $\tau(\*x_1, \*z) - \tau(\*x_0, \*z) = \3E[Y_{\*x_1} - Y_{\*x_0} \mid \*z]$, is the difference of two CAPOs, and the bounds $U_{\*z}$ and $L_{\*z}$ combine with \ref{ass:str} into its sharp identified set.
\begin{restatable}{theorem}{thmcontrast}\label{thm:contrast}
Let $\*x_1, \*x_0 \in \3X$, including treatments outside $S_{\*z}$, let $r \triangleq \rho_{\mathrm{T}}(\*x_1, \*x_0)$, and let $\1C(\*x_1, \*x_0, \*z)$ denote the interval
\begin{equation}\label{eq:contrast}
\begin{split}
    \Bigl[\, &-\3E\bigl[\min\{U_{\*z}(\*x_0, \*X') - L_{\*z}(\*x_1, \*X'),\ r\} \bigm| \*z\bigr], \\
    &\3E\bigl[\min\{U_{\*z}(\*x_1, \*X') - L_{\*z}(\*x_0, \*X'),\ r\} \bigm| \*z\bigr] \,\Bigr].
\end{split}
\end{equation}
\begin{enumerate}[label=(\alph*),leftmargin=*,topsep=2pt,itemsep=1pt]
    \item (\emph{Validity}) Under \ref{ass:bounded}--\ref{ass:smoothness_X}, $\tau(\*x_1, \*z) - \tau(\*x_0, \*z) \in \1C(\*x_1, \*x_0, \*z)$.
    \item (\emph{Sharpness}) Under \ref{ass:bounded}, \ref{ass:str}, and \ref{ass:sts}, if $\mu(\cdot, \*z)$ satisfies Eq.~\ref{eq:consistency}, then $\1C(\*x_1, \*x_0, \*z)$ is exactly the identified set of the CATE over the SCMs in $\2M(\rho_{\mathrm{T}}, \rho_{\mathrm{S}})$ that reproduce $P(\*X, \*Z, Y)$.
\end{enumerate}
\end{restatable}
Differencing the two CAPO bounds of Thm.~\ref{thm:cate_bound} gives the upper value $\3E[U_{\*z}(\*x_1, \*X') - L_{\*z}(\*x_0, \*X') \mid \*z]$, which is valid but loose, because the upper and lower bounds are attained by different ETT functions. The cap at $r$ is exactly the coupling that differencing misses: \ref{ass:str} ties the two treatment levels together within each group of units. \ref{ass:str} alone gives only the data-free interval $\pm r$ \citep[Prop.~3.2]{kim2018identification}, and \citet{kim2018identification} obtain informative CATE bounds by adding monotonicity. Thm.~\ref{thm:contrast} is informative wherever $U_{\*z}(\*x_1, \cdot) - L_{\*z}(\*x_0, \cdot) < r$ on a set of received treatments of positive probability, and needs no monotonicity.

%% file: section5.tex
\input{table_main_experiments}

\begin{figure*}[t]
\centering
\includegraphics[width=\linewidth]{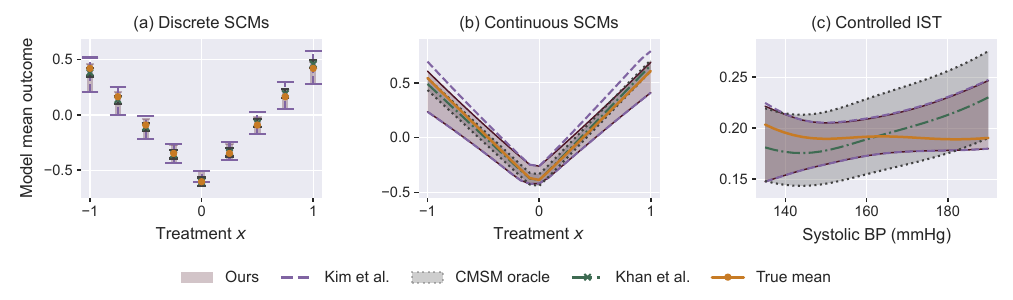}
\caption{Example bounds and true means for discrete, smooth V-shaped continuous, and controlled IST models. Wine and gray shading denote ours and oracle CMSM; dashed and dotted boundaries denote Kim and CMSM. Discrete intervals use bars. Panels have separate outcome scales; example selection is specified in App.~\ref{sec:experiment_protocol}.}
\label{fig:section5_intervals}
\end{figure*}

\section{Experiments}\label{sec:experiments}
We compare our finite-reference bound (Cor.~\ref{cor:finite}) with the combined STR+STS bound of \citet{kim2018identification}, the continuous marginal sensitivity model (CMSM; \citealp{jesson2022scalable}), and the Lipschitz-envelope adaptation of \citet{khan2024offpolicy}. Table~\ref{tab:main_experiments} reports coverage of known model means and interval width divided by the outcome range, averaged over queries within each setting. Ours and Kim receive identical regressions, Lipschitz costs, and outcome ranges. Kim uses the per-pair endpoints; ours searches the supported reference points before clipping. CMSM uses one target-based coverage oracle per SCM or IST; Khan is evaluated under violated unconfoundedness. Experimental setup and calibration details are provided in App.~\ref{sec:experiment_protocol}.

\textbf{Discrete SCMs.} We randomly draw $500$ SCMs with $U\in\{-1,+1\}$ and $K\in\{5,7,9\}$ equally spaced treatments on $[-1,1]$. Assignment follows $P(X=x_j\mid U=u)=(1+\rho ux_j
)/K$, and outcomes are affine rescalings of $b|x|+cU$ to the attained range $[-1,1]$. We perturb the response and assignment parameters. Our mean width is $19.17\%$ smaller than Kim's, with full coverage for both. The gain comes from reference points near the interior response minimum, which Kim's per-pair endpoints $\{x,x'\}$ omit.

\textbf{Continuous SCMs.} We draw $250$ SCMs from each of two families, a V-shaped response and its smooth approximation. Both use a binary hidden confounder and the density $p(x\mid u)=(1+\rho ux)/2$ on $[-1,1]$. We perturb response scales, minimum locations, latent weights, and assignment strength. Across $12500$ queries, our pooled mean width is $18.70\%$ smaller than
Kim's; both methods achieve full coverage.

\textbf{Controlled IST.} Within three age tertiles of the International Stroke Trial (\citealp{carolei1997international}), consciousness defines a hidden categorical confounder. Smoothed binary death risks define a new SCM with known, overlapping assignment on $25$ blood-pressure values per group. Target effects are computed from this model. Across $75$ targets, our pooled width reduction is $1.34\%$, with full coverage for ours and Kim. Age-group reductions are $0.47\%$, $1.32\%$, and $1.76\%$. At the first, middle, and last query of each age group, a linear program over the ETT matrix, using only the observed diagonal, STR, STS, and the outcome range, reproduces our endpoints to within $10^{-8}$.

Fig.~\ref{fig:section5_intervals} shows visible tightening over Kim in the synthetic examples and nearly coincident intervals in IST. Between oracle CMSM and ours, neither is uniformly better: CMSM is narrower in all synthetic families, while ours is narrower in all three IST groups. CMSM's full coverage, however, relies on a target-based oracle $\Lambda^\star$, calibrated using the known interventional target. Khan assumes unconfoundedness and, as expected, collapses to the observational regression.

%% file: table_main_experiments.tex
\begin{table*}[t]
\centering
\caption{Query-averaged coverage and interval width normalized by the outcome range. Reduction compares pooled mean widths with Kim. CMSM uses a per-model coverage oracle, so its coverage holds by construction.}
\label{tab:main_experiments}
\begin{tabular}{lrrrrrrrrr}
\toprule
 & \multicolumn{2}{c}{Ours} & \multicolumn{2}{c}{Kim et al.} & \multicolumn{2}{c}{CMSM oracle} & \multicolumn{2}{c}{Khan et al.} & Reduction \\
\cmidrule(lr){2-3}\cmidrule(lr){4-5}\cmidrule(lr){6-7}\cmidrule(lr){8-9}
Setting & cov. & width & cov. & width & cov. & width & cov. & width & vs. Kim (\%) \\
\midrule
Discrete SCMs & 1.00 & 0.11181 & 1.00 & 0.13832 & 1.00 & 0.04810 & 0.14 & 0.00000 & 19.17 \\
Continuous SCMs & 1.00 & 0.10627 & 1.00 & 0.13071 & 1.00 & 0.04841 & 0.04 & 0.00000 & 18.70 \\
IST & 1.00 & 0.03726 & 1.00 & 0.03777 & 1.00 & 0.07558 & 0.04 & 0.00000 & 1.34 \\
\bottomrule
\end{tabular}
\end{table*}

%% file: section6.tex
\vspace{-0.05in}
\section{Conclusions}
We studied partial identification of CAPOs and CATEs for continuous treatments under unobserved confounding. Lipschitz response functions and Lipschitz propensity scores given the latent confounder imply the STR and STS conditions of \citet{kim2018identification}, with a selection cost that depends only on the observed treatment density. Under STR and STS with general costs, bounding through every diagonal reference point gives the sharp identified set of the CAPO at every treatment, inside or outside the observational support, whereas the two-point combined bound of \citet{kim2018identification} is not sharp in general. Capping the difference of the upper and lower bounds at the response cost gives the sharp identified set of the CATE, without any monotonicity assumption.

%% file: appendix_e.tex
\section{Related Work} \label{appendix:related}

This section expands the discussion of the Introduction, placing our bound among the main lines of work on identification, partial identification, and sensitivity analysis, and stating precisely how it relates to each.

\textbf{Identification and adjustment.} When the target effect is point-identified, for example under the back-door or ignorability condition \citep{pearl:2k, rosenbaum:rub83}, it can be recovered from observational data by covariate adjustment \citep{rosenbaum:rub83, bang2005doubly} or, in the sequential and off-policy settings, by importance-weighted and doubly robust estimators \citep{dudik2011doubly, li2015toward, munos2016safe, thomas2016data}. Whether a given effect is identified at all is itself settled by the do-calculus, which is complete for interventional distributions \citep{huang:val06, shpitser:pea06a, bareinboim2022pearl}. Our setting begins where these tools end: unobserved confounding leaves the CAPO non-identified, so no adjustment functional exists, and the effect can only be bounded.

\textbf{Partial identification for discrete treatments.} The classical response to non-identification is to bound the effect over all models compatible with the data \citep{manski1990nonparametric, balke:pea97, frangakis:rub02}. \citet{manski1990nonparametric} gives the worst-case interval that assumes nothing about the unobserved response, which our Cor.~\ref{cor:manski} recovers as the no-smoothness corner of the family. \citet{balke:pea97} tighten it under an instrument by linear programming over a discretized latent domain, and \citet{zhang2022partial} extend such programs to counterfactual queries and to Monte-Carlo bounds for discrete outcomes; \citet{manski1998monotone} instead impose a monotone instrumental-variable structure. These bounds are sharp, but they rest on enumerating a finite latent domain, and that enumeration does not survive the move to continuous treatments and outcomes, which is the regime we target.

\textbf{Sensitivity analysis and bounded confounding.} A complementary line fixes a quantitative limit on how far confounding can distort the treatment mechanism and reports the resulting interval. For binary and discrete treatments, the marginal sensitivity model and its relatives bound the ratio between the observed and the latent-conditional propensity \citep{kallus2018confounding, namkoong2020off}. For continuous treatments, \citet{jesson2022scalable} introduce the continuous marginal sensitivity model (CMSM), which we adopt as our main baseline, and bound both the CAPO $\3E[Y_{\*x} \mid \*z]$ and its average. \citet{frauen2023sharp} derive sharp bounds, pointwise in the covariates, for a generalized sensitivity model that contains the CMSM, and \citet{baitairian2025sharp} give sharp closed-form CMSM bounds with a doubly robust estimator. \citet{marmarelis2023partial} instead bound how fast the log density ratio may change with the treatment (the $\delta$-MSM), \citet{dalal2025partial} let the sensitivity parameter vary with the exposure level, and \citet{bonvini2022sensitivity} develop propensity-based, outcome-based, and subset-confounding sensitivity models for marginal structural models. Several of these methods bound conditional, pointwise effects; what separates them from ours is the form of the restriction. They constrain the treatment density relative to its observed value and are stated where that density is positive, so their intervals are built from the observed outcomes at the query and do not extend to treatments outside the support. Our restriction is instead on the smoothness of the mechanisms.

\textbf{Identification under parametric structure.} When the mechanisms are constrained, the effect can sometimes be point-identified despite confounding. Linear structural equation models identify effects through path coefficients \citep{wright1934method} and auxiliary variables \citep{chen2017identification}, and the deconfounder exploits many causes sharing a latent factor \citep{wang2019blessings}. We make no such parametric commitment. Our inputs are smoothness costs, and \citet{zhang2026causal} show that smooth-mechanism SCMs are dense among the confounded models compatible with the data, matching the observational and every interventional law to arbitrary accuracy. Smoothness as such therefore rules out no observational law; the substantive, and untestable beyond Cor.~\ref{cor:test}, assumption is the size of the constants, which we treat as sensitivity parameters.

\textbf{Optimization-based and generative bounds.} A further line computes bounds numerically, by optimizing over a class of structural models consistent with the data. For continuous treatments with an instrument, \citet{kilbertus2020class} optimize over response functions spanned by a chosen basis, \citet{padh2023stochastic} extend this to instruments and leaky mediators by stochastic constrained optimization, and \citet{gunsilius2020path} solves the resulting infinite-dimensional programs by path sampling. \citet{balazadeh2022partial} fit implicit generative models to general causal graphs and bound average treatment derivatives by optimizing over them, and neural or generative SCMs are used in the same way elsewhere \citep{kocaoglu2018causalgan, xia2021causal, xia2022neural, nasr2023counterfactual}. \citet{tan2024consistency} analyze when such neural estimates are consistent, and make explicit that an informative bound over continuous domains is impossible without positing a smoothness assumption, since two models can match the data yet disagree arbitrarily on a rarely observed treatment. The bounds in this line are valid relative to the chosen function class and are the output of a nonconvex optimization; ours is an explicit bound under stated Lipschitz constants.

\textbf{Smoothness assumptions on counterfactual means.} Closest to our work is a line in econometrics that bounds counterfactual means by limiting their rate of change. \citet{hausman2016individual} bound each individual's income effect to partially identify average welfare under unrestricted heterogeneity. Following them, \citet{kim2018identification} study $g(x, x') = \3E[Y_x \mid X = x']$ for a scalar treatment, assume that $g$ is Lipschitz in the intervened treatment (smooth treatment response, STR: $|g(x, x') - g(\tilde x, x')| \le b|x - \tilde x|$) or in the received treatment (smooth treatment selection, STS: $|g(x, x'_0) - g(x, x'_1)| \le a|x'_0 - x'_1|$), and derive closed-form bounds on $\3E[Y_x]$ under each assumption, under both, and in combination with monotonicity; their bounds under monotone treatment response and selection (Props.~4.1--4.2) already route through intermediate treatment levels, with a monotone rather than Lipschitz selection step. Our ETT $g_{\*z}(\*x, \*x')$ is their $g$ with covariates $\*z$, our primary assumptions \ref{ass:str} and \ref{ass:sts} are theirs with general costs, and our STR and STS bounds correspond to their propositions:
\begin{center}
\begin{tabular}{ll}
\toprule
This paper & \citet{kim2018identification} \\
\midrule
Lem.~\ref{lem:side_bound} (intervened treatment) & Prop.~3.1 (STR) \\
Lem.~\ref{lem:cond_bound} (received treatment) & Prop.~D.1 (STS) \\
Cor.~\ref{cor:finite} with $\3I = \{\*x, \*x'\}$ at each $\*x'$ & Prop.~D.3 (STR and STS) \\
Thm.~\ref{thm:sharp} (all reference points) & -- \\
\bottomrule
\end{tabular}
\end{center}
Our contribution is therefore not a new bounding principle: the STR and STS bounds and their combination go back to \citet{kim2018identification}, and the combined bound is a standard Lipschitz-extension bound \citep{mcshane1934extension}. Relative to their work, five things are new. \emph{First, sharpness.} Prop.~D.3 of \citet{kim2018identification} is stated as sharp, but its proof constructs a data-generating process with a constant regression, for which intermediate reference points never help. For a nonconstant regression the two-point bound is not sharp (Prop.~\ref{prop:kimloose}), whereas our bound is the sharp identified set (Thm.~\ref{thm:sharp}); we also give the exact condition under which the full reference set is strictly tighter (Cor.~\ref{cor:strict}) and the testable implication of the model (Cor.~\ref{cor:test}). \emph{Second, the selection condition.} \citet{kim2018identification} assume STS directly and motivate it informally; we derive it (Prop.~\ref{prop:structural}) from a Lipschitz assignment density given the latent confounder, via a total-variation bound on the latent posteriors (Lem.~\ref{lem:posterior_smooth}), with the identified, density-adaptive cost $\kappa^{\ast}_{\*z}$, which needs no overlap floor. \emph{Third, the scope.} We allow covariates, with effects conditional on them, and vector-valued treatments, and we use the bounded outcome range inside the bound, whereas their main results concern a scalar treatment without covariates and with unbounded outcomes. \emph{Fourth, extrapolation.} Their STR bound is valid at treatments outside the support, but their combined bound then reduces to it, because the STS term needs the regression at the query; our bound still uses both steps there and remains the sharp identified set (Thm.~\ref{thm:sharp}). \emph{Fifth, the CATE.} Under STR alone their CATE bound is the data-free interval $\pm b|t - t'|$ (Prop.~3.2), and their informative CATE bounds add monotonicity; under STR and STS, Thm.~\ref{thm:contrast} gives the sharp identified set of the CATE without it.

\textbf{Smoothness and off-policy evaluation beyond overlap.} \citet{khan2024offpolicy} bound the off-policy value under a Lipschitz smoothness assumption through closed-form linear programs, targeting the unconfounded setting with limited or no overlap between the logging and target policies. Our setting differs in that confounding renders the effect non-identified even under full overlap. In the unconfounded limit our STR bound reduces to their bound (Cor.~\ref{cor:khan}), so their result is the no-confounding case of ours; the STS bound, absent from their analysis, is what accounts for confounding breaking the equality of the observed regression and the true effect.

%% file: appendix_a.tex
\section{Proofs} \label{appendix:a}
This appendix collects the proofs of the theoretical results in the main text. App.~\ref{app:proofs_structural} proves the structural sufficient conditions (Lem.~\ref{lem:posterior_smooth} and Prop.~\ref{prop:structural}) that connect \ref{ass:smoothness_Y} and \ref{ass:smoothness_X} to the STR and STS conditions; App.~\ref{app:proofs_side} proves the STR and STS bounds (Lems.~\ref{lem:side_bound} and~\ref{lem:cond_bound}) and their unconfounded special case (Cor.~\ref{cor:khan}); App.~\ref{app:proofs_combined} proves the combined bound (Lem.~\ref{lem:rect_anchor} and Thm.~\ref{thm:cate_bound}) together with its finite-reference corollary (Cor.~\ref{cor:finite}) and the strictness and Manski corollaries (Cors.~\ref{cor:strict} and~\ref{cor:manski}); and App.~\ref{app:proofs_sharp} proves sharpness (Thms.~\ref{thm:sharp} and~\ref{thm:contrast}, Cor.~\ref{cor:test}) and the comparison with the bound of \citet{kim2018identification} (Props.~\ref{prop:kim} and~\ref{prop:kimloose}). All validity proofs follow the same playbook: the ETT function $g_{\*z}(\*x, \*x')$ is bounded pointwise by chaining consistency on the diagonal with the STR and STS inequalities, and the pointwise bound is then averaged over the received treatment $\*X' \sim P(\*X \mid \*z)$ through the tower property $\tau(\*x, \*z) = \3E[g_{\*z}(\*x, \*X') \mid \*z]$ of Eq.~\ref{eq:cate_expand}. The sharpness proofs run in the reverse direction: they exhibit an admissible ETT function that attains the bound and realize it by an SCM that reproduces the observational law.

Throughout, statements about a fixed covariate value $\*z$ hold for $P$-almost every $\*z$, with regular versions of the conditional laws, and all integrated functions of $(\*z, \*x, \*x')$ are jointly measurable. Apart from Lem.~\ref{lem:posterior_smooth} and Prop.~\ref{prop:structural}, the proofs use only \ref{ass:bounded}, \ref{ass:str}, \ref{ass:sts}, and consistency, with arbitrary costs; by Prop.~\ref{prop:structural} and the convention of Sec.~\ref{sec:structural}, this covers \ref{ass:bounded}--\ref{ass:smoothness_X}, with $S_{\*z}$ replaced by $S^{+}_{\*z}$, which carries all of $P(\*X \mid \*z)$ and whose points have neighbourhoods of positive mass (see the proof of Prop.~\ref{prop:structural}).

\subsection{Structural sufficient conditions}\label{app:proofs_structural}
\lemposteriorsmooth*
\begin{proof}
Fix $\*z$ and write
\[
q \triangleq P(\*U \mid \*z), \qquad
p_{\*u}(\*x) \triangleq p(\*x \mid \*z, \*u), \qquad
p(\*x) \triangleq p(\*x \mid \*z) = \int p_{\*u}(\*x)\, dq(\*u), \qquad
\Delta \triangleq \|\*x_0' - \*x_1'\|.
\]
Each $p_{\*u}$ is $L_X$-Lipschitz by \ref{ass:smoothness_X}. Since $p$ is an average of these functions, it is $L_X$-Lipschitz as well:
\[
|p(\*x_0') - p(\*x_1')|
\le
\int |p_{\*u}(\*x_0') - p_{\*u}(\*x_1')|\, dq(\*u)
\le
L_X \Delta.
\]
By Bayes' rule, whenever $p(\*x') > 0$,
\[
dP(\*u \mid \*x', \*z) = \frac{p_{\*u}(\*x')}{p(\*x')}\, dq(\*u).
\]
Hence the total variation between the two latent posteriors is
\[
2\,\mathrm{TV}_{\*u}(\*x_0', \*x_1'; \*z)
=
\int \Bigl|\frac{p_{\*u}(\*x_0')}{p(\*x_0')} - \frac{p_{\*u}(\*x_1')}{p(\*x_1')}\Bigr|\, dq(\*u).
\]
Adding and subtracting $p_{\*u}(\*x_1') / p(\*x_0')$ inside the absolute value and applying the triangle inequality, we obtain
\begin{align*}
    2\,\mathrm{TV}_{\*u}(\*x_0', \*x_1'; \*z)
    &\le \frac{1}{p(\*x_0')}\int |p_{\*u}(\*x_0') - p_{\*u}(\*x_1')|\, dq(\*u) + \Bigl|\frac{1}{p(\*x_0')} - \frac{1}{p(\*x_1')}\Bigr| \int p_{\*u}(\*x_1')\, dq(\*u).
\end{align*}
By \ref{ass:smoothness_X}, the first term is at most $L_X \Delta / p(\*x_0')$. Since $\int p_{\*u}(\*x_1')\, dq(\*u) = p(\*x_1')$, the second term equals
\[
\Bigl|\frac{1}{p(\*x_0')} - \frac{1}{p(\*x_1')}\Bigr|\, p(\*x_1')
=
\frac{|p(\*x_1') - p(\*x_0')|}{p(\*x_0')}
\le
\frac{L_X \Delta}{p(\*x_0')},
\]
where the inequality uses the Lipschitz property of $p$. Combining the two terms gives
\[
2\,\mathrm{TV}_{\*u}(\*x_0', \*x_1'; \*z)
\le
\frac{L_X \Delta}{p(\*x_0')} + \frac{L_X \Delta}{p(\*x_0')}
=
\frac{2 L_X \Delta}{p(\*x_0')}.
\]
Exchanging the roles of $\*x_0'$ and $\*x_1'$ gives the same bound with $p(\*x_1')$ in the denominator. Taking the better of the two bounds gives Eq.~\ref{eq:posterior_smooth}. This completes the proof.
\end{proof}

\propstructural*
\begin{proof}
We first prove \ref{ass:str}. Fix $\*x, \tilde{\*x} \in \3X$ and $\*x' \in S^{+}_{\*z}$. By Eq.~\ref{eq:ett}, the two ETT values integrate different response functions against the same latent posterior:
\[
g_{\*z}(\*x, \*x') - g_{\*z}(\tilde{\*x}, \*x')
=
\int \bigl(h_{\*u}(\*x, \*z) - h_{\*u}(\tilde{\*x}, \*z)\bigr)\, dP(\*u \mid \*x', \*z).
\]
By \ref{ass:smoothness_Y}, $|h_{\*u}(\*x, \*z) - h_{\*u}(\tilde{\*x}, \*z)| \le L_Y\|\*x - \tilde{\*x}\|$ for every $\*u$. Since the bound holds for every stratum, it holds for the average:
\[
|g_{\*z}(\*x, \*x') - g_{\*z}(\tilde{\*x}, \*x')|
\le
L_Y\|\*x - \tilde{\*x}\|
=
\rho_{\mathrm{T}}(\*x, \tilde{\*x}).
\]
This is \ref{ass:str}.

We next prove \ref{ass:sts}. Fix $\*x \in \3X$ and $\*x_0', \*x_1' \in S^{+}_{\*z}$. By Eq.~\ref{eq:ett}, the two ETT values share the response functions and differ only through the signed measure
\[
d\Delta(\*u) \triangleq dP(\*u \mid \*x_0', \*z) - dP(\*u \mid \*x_1', \*z),
\]
which has zero total mass. Let $c$ be the midpoint of $\3Y$. By \ref{ass:bounded}, each $h_{\*u}(\*x, \*z) \in \3Y$, so $|h_{\*u}(\*x, \*z) - c| \le \tfrac{1}{2}\operatorname{diam}(\3Y)$. Since $\Delta$ has zero total mass, subtracting $c$ from the integrand does not change the integral, and therefore
\[
\bigl|g_{\*z}(\*x, \*x_0') - g_{\*z}(\*x, \*x_1')\bigr|
=
\Bigl|\int \bigl(h_{\*u}(\*x, \*z) - c\bigr)\, d\Delta(\*u)\Bigr|
\le
\operatorname{diam}(\3Y)\,\mathrm{TV}_{\*u}(\*x_0', \*x_1'; \*z).
\]
Since $p(\*x_0' \mid \*z) > 0$ and $p(\*x_1' \mid \*z) > 0$ on $S^{+}_{\*z}$, Lem.~\ref{lem:posterior_smooth} bounds the right side by $\kappa_{\*z}(\*x_0', \*x_1')$; this is Eq.~\ref{eq:transport_xp}. To pass to the chain closure, let $\*a = \*a_0, \dots, \*a_k = \*b$ be a chain in $S^{+}_{\*z}$. Applying Eq.~\ref{eq:transport_xp} to each link and the triangle inequality, we obtain
\[
|g_{\*z}(\*x, \*a) - g_{\*z}(\*x, \*b)|
\le
\sum_{i} |g_{\*z}(\*x, \*a_i) - g_{\*z}(\*x, \*a_{i+1})|
\le
\sum_{i} \kappa_{\*z}(\*a_i, \*a_{i+1}).
\]
Taking the infimum over chains gives $|g_{\*z}(\*x, \*a) - g_{\*z}(\*x, \*b)| \le \kappa^{\ast}_{\*z}(\*a, \*b)$, which is \ref{ass:sts} with $\kappa^{\ast}_{\*z}$.

We next verify that the cost is regular. The cost $\kappa^{\ast}_{\*z}$ is symmetric and vanishes on the diagonal. Concatenating a chain from $\*a$ to $\*b$ with a chain from $\*b$ to $\*c$ gives a chain from $\*a$ to $\*c$, so
\[
\kappa^{\ast}_{\*z}(\*a, \*c) \le \kappa^{\ast}_{\*z}(\*a, \*b) + \kappa^{\ast}_{\*z}(\*b, \*c).
\]
Hence $\kappa^{\ast}_{\*z}$ is a pseudometric. It is continuous on $S^{+}_{\*z}$. Indeed, for $\*a \in S^{+}_{\*z}$ and $\*b \to \*a$, the one-link chain gives
\[
\kappa^{\ast}_{\*z}(\*a, \*b)
\le
\kappa_{\*z}(\*a, \*b)
\le
\frac{L_X\operatorname{diam}(\3Y)\,\|\*a - \*b\|}{p(\*a \mid \*z)}
\to 0,
\]
and the triangle inequality gives
\[
|\kappa^{\ast}_{\*z}(\*a, \*b) - \kappa^{\ast}_{\*z}(\*a', \*b')|
\le
\kappa^{\ast}_{\*z}(\*a, \*a') + \kappa^{\ast}_{\*z}(\*b, \*b').
\]

Finally, we show that $S^{+}_{\*z}$ has full measure. Since $P(\*X \mid \*z)$ has density $p(\cdot \mid \*z)$,
\[
P\bigl(p(\*X \mid \*z) = 0 \bigm| \*z\bigr)
=
\int_{\{p = 0\}} p(\*x \mid \*z)\, d\*x
=
0,
\]
so $S^{+}_{\*z}$ carries all of $P(\*X \mid \*z)$. The proofs of Thms.~\ref{thm:cate_bound} and~\ref{thm:sharp} use the support only through this property and through the fact that each of its points has neighbourhoods of positive mass, which $S^{+}_{\*z} \subseteq S_{\*z}$ inherits. Hence both theorems hold with $S^{+}_{\*z}$ in place of $S_{\*z}$, and no integral changes. This completes the proof.
\end{proof}

\subsection{STR and STS bounds}\label{app:proofs_side}
\lemsidebound*
\begin{proof}
By the tower property over the received treatment $\*X' \sim P(\cdot \mid \*z)$,
\[
\3E[Y_{\*x} \mid \*z] = \3E\bigl[g_{\*z}(\*x, \*X') \bigm| \*z\bigr].
\]
It therefore suffices to bound $g_{\*z}(\*x, \*x')$ pointwise for $\*x' \in S_{\*z}$. By consistency,
\[
g_{\*z}(\*x', \*x') = \mu(\*x', \*z), \qquad \*x' \in S_{\*z}.
\]
Applying \ref{ass:str} at the received treatment $\*x'$, between the intervened treatments $\*x$ and $\*x'$, we obtain
\[
|g_{\*z}(\*x, \*x') - \mu(\*x', \*z)|
=
|g_{\*z}(\*x, \*x') - g_{\*z}(\*x', \*x')|
\le
\rho_{\mathrm{T}}(\*x, \*x').
\]
Moreover, $g_{\*z}(\*x, \*x') \in [\underline{y}, \overline{y}]$ by \ref{ass:bounded}, being an average of values of $f_Y$. Combining the two constraints, for every $\*x' \in S_{\*z}$,
\[
\bigl(\mu(\*x', \*z) - \rho_{\mathrm{T}}(\*x, \*x')\bigr) \vee \underline{y}
\;\le\;
g_{\*z}(\*x, \*x')
\;\le\;
\bigl(\mu(\*x', \*z) + \rho_{\mathrm{T}}(\*x, \*x')\bigr) \wedge \overline{y}.
\]
Substituting $\*x' = \*X'$ and taking expectations given $\*z$, monotonicity of the expectation yields
\[
\underline{w}(\*x, \*z)
\le
\3E\bigl[g_{\*z}(\*x, \*X') \bigm| \*z\bigr]
=
\3E[Y_{\*x} \mid \*z]
\le
\overline{w}(\*x, \*z),
\]
where the two sides are the definitions Eq.~\ref{eq:side_halfwidth}. This proves Eq.~\ref{eq:side_anchor}.
\end{proof}

\begin{corollary}\label{cor:khan}
Under \ref{ass:bounded}--\ref{ass:smoothness_X}, suppose additionally there is no unobserved confounding, $Y_{\*x} \perp \*X \mid \*Z$. Then $\mu(\*x, \*z) = \tau(\*x, \*z)$ on $S_{\*z}$ and the ETT is constant in the received treatment, $\3E[Y_{\*x} \mid \*x', \*z] = \tau(\*x, \*z)$ for every $\*x' \in S_{\*z}$. Every observed treatment therefore constrains the \emph{same} constant $\tau(\*x, \*z)$, so the average over $P(\*X \mid \*z)$ that defines the bounds $\underline{w}(\*x, \*z), \overline{w}(\*x, \*z)$ can be replaced by the intersection of the per-treatment intervals, sharpening it to the Lipschitz bound of \citet{khan2024offpolicy}
\begin{equation}\label{eq:khan_envelope}
\begin{split}
    \sup_{\*x' \in S_{\*z}} \bigl[\mu(\*x', \*z) - \rho_{\mathrm{T}}(\*x, \*x')\bigr] \le \tau(\*x, \*z) \\
    \le \inf_{\*x' \in S_{\*z}} \bigl[\mu(\*x', \*z) + \rho_{\mathrm{T}}(\*x, \*x')\bigr].
\end{split}
\end{equation}
\end{corollary}
\begin{proof}
Under $Y_{\*x} \perp \*X \mid \*Z$, conditioning on the received treatment does not change the mean of the potential outcome, so for every $\*x' \in S_{\*z}$,
\[
g_{\*z}(\*x, \*x') = \3E[Y_{\*x} \mid \*x', \*z] = \3E[Y_{\*x} \mid \*z] = \tau(\*x, \*z),
\]
independent of $\*x'$. Taking $\*x' = \*x \in S_{\*z}$ and using consistency $g_{\*z}(\*x, \*x) = \mu(\*x, \*z)$, we obtain $\mu(\*x, \*z) = \tau(\*x, \*z)$.

The pointwise bound in the proof of Lem.~\ref{lem:side_bound} now constrains this constant: for every observed $\*x' \in S_{\*z}$,
\[
\tau(\*x, \*z) \in \bigl[\mu(\*x', \*z) - \rho_{\mathrm{T}}(\*x, \*x'),\ \mu(\*x', \*z) + \rho_{\mathrm{T}}(\*x, \*x')\bigr] \cap \3Y.
\]
Since every $\*x'$ constrains the same constant $\tau(\*x, \*z)$, the average over $P(\*X \mid \*z)$ may be replaced by the intersection over $\*x' \in S_{\*z}$. Taking the supremum of the lower ends and the infimum of the upper ends gives Eq.~\ref{eq:khan_envelope}.

It remains to show that each side is attained at a single observed treatment. Chaining the upper bound through an intermediate observed treatment $\*m \in S_{\*z}$, first from $\*x'$ to $\*m$ and then from $\*m$ to $\*x$, gives the upper value $\mu(\*x', \*z) + \rho_{\mathrm{T}}(\*m, \*x') + \rho_{\mathrm{T}}(\*x, \*m)$. By the triangle inequality for $\rho_{\mathrm{T}}$,
\[
\mu(\*x', \*z) + \rho_{\mathrm{T}}(\*m, \*x') + \rho_{\mathrm{T}}(\*x, \*m)
\ge
\mu(\*x', \*z) + \rho_{\mathrm{T}}(\*x, \*x'),
\]
so passing through $\*m$ only relaxes the bound, and the nearest effective observed treatment determines it (the no-interaction property). The lower side is symmetric. This completes the proof.
\end{proof}
\lemcondbound*
\begin{proof}
The proof follows that of Lem.~\ref{lem:side_bound}, with \ref{ass:sts} in place of \ref{ass:str}. By the tower property,
\[
\3E[Y_{\*x} \mid \*z] = \3E\bigl[g_{\*z}(\*x, \*X') \bigm| \*z\bigr],
\]
while $\mu(\*x, \*z)$ is constant in $\*X'$. Since $\*x \in S_{\*z}$, consistency gives
\[
g_{\*z}(\*x, \*x) = \mu(\*x, \*z).
\]
Applying \ref{ass:sts} at the intervened treatment $\*x$, between the received treatments $\*x'$ and $\*x$, we obtain
\[
|g_{\*z}(\*x, \*x') - \mu(\*x, \*z)|
=
|g_{\*z}(\*x, \*x') - g_{\*z}(\*x, \*x)|
\le
\rho_{\mathrm{S}}(\*x', \*x).
\]
As before, $g_{\*z}(\*x, \*x') \in [\underline{y}, \overline{y}]$. Hence, for every $\*x' \in S_{\*z}$,
\[
\bigl(\mu(\*x, \*z) - \rho_{\mathrm{S}}(\*x', \*x)\bigr) \vee \underline{y}
\;\le\;
g_{\*z}(\*x, \*x')
\;\le\;
\bigl(\mu(\*x, \*z) + \rho_{\mathrm{S}}(\*x', \*x)\bigr) \wedge \overline{y}.
\]
Substituting $\*x' = \*X'$ and using monotonicity of the expectation gives Eq.~\ref{eq:cond_side_anchor}, the two sides being Eq.~\ref{eq:cond_halfwidth}. This completes the proof.
\end{proof}

\subsection{Combined bounds}\label{app:proofs_combined}
\begin{lemma}\label{lem:product}
Let $\rho_{\mathrm{T}}$ and $\rho_{\mathrm{S}}$ be pseudometrics on $\3X$ and $S_{\*z}$, and $d_{\*z}((\*x, \*x'), (\*a, \*b)) \triangleq \rho_{\mathrm{T}}(\*x, \*a) + \rho_{\mathrm{S}}(\*x', \*b)$. A function $g : \3X \times S_{\*z} \to \3R$ satisfies the inequalities of \ref{ass:str} and \ref{ass:sts} with these costs if and only if $|g(\*p) - g(\*q)| \le d_{\*z}(\*p, \*q)$ for all $\*p, \*q \in \3X \times S_{\*z}$.
\end{lemma}
\begin{proof}
We first show that a $d_{\*z}$-Lipschitz function satisfies the two conditions. Let $g$ be $d_{\*z}$-Lipschitz. Taking $\*p = (\*x, \*x')$ and $\*q = (\tilde{\*x}, \*x')$, and using $\rho_{\mathrm{S}}(\*x', \*x') = 0$,
\[
|g(\*x, \*x') - g(\tilde{\*x}, \*x')|
\le
\rho_{\mathrm{T}}(\*x, \tilde{\*x}) + \rho_{\mathrm{S}}(\*x', \*x')
=
\rho_{\mathrm{T}}(\*x, \tilde{\*x}),
\]
which is the inequality of \ref{ass:str}. Taking $\*p = (\*x, \*x_0')$ and $\*q = (\*x, \*x_1')$, and using $\rho_{\mathrm{T}}(\*x, \*x) = 0$,
\[
|g(\*x, \*x_0') - g(\*x, \*x_1')|
\le
\rho_{\mathrm{T}}(\*x, \*x) + \rho_{\mathrm{S}}(\*x_0', \*x_1')
=
\rho_{\mathrm{S}}(\*x_0', \*x_1'),
\]
which is the inequality of \ref{ass:sts}.

We next prove the converse. Suppose $g$ satisfies the inequalities of \ref{ass:str} and \ref{ass:sts}, and let $\*p = (\*x, \*x')$ and $\*q = (\*a, \*b)$. Passing through the intermediate point $(\*a, \*x')$ and applying \ref{ass:str} to the first step and \ref{ass:sts} to the second,
\begin{equation*}
    |g(\*x, \*x') - g(\*a, \*b)| \le |g(\*x, \*x') - g(\*a, \*x')| + |g(\*a, \*x') - g(\*a, \*b)| \le \rho_{\mathrm{T}}(\*x, \*a) + \rho_{\mathrm{S}}(\*x', \*b).
\end{equation*}
The right side is $d_{\*z}(\*p, \*q)$, so $g$ is $d_{\*z}$-Lipschitz. This completes the proof.
\end{proof}

\lemrectanchor*
\begin{proof}
By Prop.~\ref{prop:structural} and the convention of Sec.~\ref{sec:structural}, the ETT function satisfies \ref{ass:str} and \ref{ass:sts} with the costs of Prop.~\ref{prop:structural}. By Lem.~\ref{lem:product}, it is therefore $d_{\*z}$-Lipschitz:
\[
|g_{\*z}(\*p) - g_{\*z}(\*q)| \le d_{\*z}(\*p, \*q), \qquad \*p, \*q \in \3X \times S_{\*z}.
\]
Now take $\*p = (\*x, \*x')$ and $\*q = (\*s, \*s)$ with $\*s \in S_{\*z}$. By consistency, $g_{\*z}(\*s, \*s) = \mu(\*s, \*z)$, and by the definition of the product cost,
\[
d_{\*z}\bigl((\*x, \*x'), (\*s, \*s)\bigr) = \rho_{\mathrm{T}}(\*x, \*s) + \rho_{\mathrm{S}}(\*x', \*s) = \rho(\*s, \*x').
\]
Substituting both into the Lipschitz inequality gives
\[
\bigl|\3E[Y_{\*x} \mid \*x', \*z] - \mu(\*s, \*z)\bigr|
=
|g_{\*z}(\*x, \*x') - g_{\*z}(\*s, \*s)|
\le
\rho(\*s, \*x'),
\]
which is Eq.~\ref{eq:rect_anchor}. This completes the proof.
\end{proof}

\paragraph{Optimality of the L-shaped path.} Lem.~\ref{lem:rect_anchor} reaches the reference point along a two-step, L-shaped path: one step in the intervened treatment and one in the received treatment (Fig.~\ref{fig:noint_L}, blue). Chaining the two conditions along any other path that moves one coordinate at a time, such as the zig-zag of Fig.~\ref{fig:noint_L} (grey), cannot give a smaller cost. To see this, let
\[
(\*x, \*x') = (\*a_0, \*b_0), (\*a_1, \*b_1), \dots, (\*a_k, \*b_k) = (\*s, \*s)
\]
be such a path in $\3X \times S_{\*z}$, where each step changes either the intervened treatment, at cost $\rho_{\mathrm{T}}(\*a_i, \*a_{i+1})$ by \ref{ass:str}, or the received treatment, at cost $\rho_{\mathrm{S}}(\*b_i, \*b_{i+1})$ by \ref{ass:sts}. Summing the steps bounds $|g_{\*z}(\*x, \*x') - \mu(\*s, \*z)|$ by
\begin{equation*}
    \sum_{i \in \mathcal{T}} \rho_{\mathrm{T}}(\*a_i, \*a_{i+1}) + \sum_{i \in \mathcal{S}} \rho_{\mathrm{S}}(\*b_i, \*b_{i+1}),
\end{equation*}
where $\mathcal{T}$ and $\mathcal{S}$ index the STR and the STS steps. The STR steps form a chain from $\*x$ to $\*s$ in the first argument, and the STS steps a chain from $\*x'$ to $\*s$ in the second, because each argument is unchanged by the steps of the other kind. By the triangle inequality for each cost,
\[
\sum_{i \in \mathcal{T}} \rho_{\mathrm{T}}(\*a_i, \*a_{i+1}) \ge \rho_{\mathrm{T}}(\*x, \*s),
\qquad
\sum_{i \in \mathcal{S}} \rho_{\mathrm{S}}(\*b_i, \*b_{i+1}) \ge \rho_{\mathrm{S}}(\*x', \*s),
\]
and the sum of the two right sides is the cost $\rho(\*s, \*x')$ of the L-shaped path. The two orders of the L, through $(\*s, \*x')$ or through $(\*x, \*s)$, cost the same. A path through several diagonal points does no better either: bounding $g_{\*z}(\*x, \*x')$ by $\mu(\*s_1, \*z)$ and then $\mu(\*s_1, \*z)$ by $\mu(\*s_2, \*z)$ costs
\[
d_{\*z}\bigl((\*x, \*x'), (\*s_1, \*s_1)\bigr) + d_{\*z}\bigl((\*s_1, \*s_1), (\*s_2, \*s_2)\bigr)
\ge
d_{\*z}\bigl((\*x, \*x'), (\*s_2, \*s_2)\bigr)
=
\rho(\*s_2, \*x'),
\]
and the bounds of Eq.~\ref{eq:envelopes} already optimize over the endpoint. So the L-shaped path to the best reference point is the cheapest way to combine the two conditions, and Thm.~\ref{thm:sharp} shows more: no argument at all can improve the resulting bound, because it is attained. The geometry is standard Lipschitz geometry; what is not standard is which bound is sharp, and that intermediate reference points are needed to reach it.

\thmcatebound*
\begin{proof}
Let $g_{\*z}$ be the ETT function of an SCM in $\2M(\rho_{\mathrm{T}}, \rho_{\mathrm{S}})$ that reproduces $P(\*X, \*Z, Y)$. Its diagonal is $\mu(\cdot, \*z)$ by consistency, and it takes values in $\3Y$ by \ref{ass:bounded}. Fix $(\*x, \*x') \in \3X \times S_{\*z}$ and a reference point $\*s \in S_{\*z}$. By Lem.~\ref{lem:product}, $g_{\*z}$ is $d_{\*z}$-Lipschitz, so consistency gives
\[
g_{\*z}(\*x, \*x')
\le
g_{\*z}(\*s, \*s) + d_{\*z}\bigl((\*x, \*x'), (\*s, \*s)\bigr)
=
\mu(\*s, \*z) + \rho(\*s, \*x').
\]
Since this holds for every $\*s \in S_{\*z}$, taking the infimum over $\*s$ and using $g_{\*z} \le \overline{y}$ gives
\[
g_{\*z}(\*x, \*x')
\le
\Bigl(\inf_{\*s \in S_{\*z}} \bigl[\mu(\*s, \*z) + \rho(\*s, \*x')\bigr]\Bigr) \wedge \overline{y}
=
U_{\*z}(\*x, \*x').
\]
The lower bound is symmetric: $g_{\*z}(\*x, \*x') \ge \mu(\*s, \*z) - \rho(\*s, \*x')$ for every $\*s \in S_{\*z}$, and $g_{\*z} \ge \underline{y}$, so
\[
g_{\*z}(\*x, \*x')
\ge
\Bigl(\sup_{\*s \in S_{\*z}} \bigl[\mu(\*s, \*z) - \rho(\*s, \*x')\bigr]\Bigr) \vee \underline{y}
=
L_{\*z}(\*x, \*x').
\]
This proves $L_{\*z} \le g_{\*z} \le U_{\*z}$ on $\3X \times S_{\*z}$.

We next average over the received treatment. The bounds $U_{\*z}$ and $L_{\*z}$ are bounded, so the expectations below exist. Substituting $\*x' = \*X' \sim P(\*X \mid \*z)$ in the pointwise bounds and using $\tau(\*x, \*z) = \3E[g_{\*z}(\*x, \*X') \mid \*z]$, monotonicity of the expectation gives
\[
\3E\bigl[L_{\*z}(\*x, \*X') \mid \*z\bigr]
\le
\3E\bigl[g_{\*z}(\*x, \*X') \mid \*z\bigr]
=
\tau(\*x, \*z)
\le
\3E\bigl[U_{\*z}(\*x, \*X') \mid \*z\bigr],
\]
which is Eq.~\ref{eq:cate_valid}. This completes the proof.
\end{proof}

\corfinite*
\begin{proof}
Since $\3I \subseteq S_{\*z}$, the infimum over $\3I$ is at least the infimum over $S_{\*z}$: for every $(\*x, \*x')$,
\[
\inf_{\*s \in S_{\*z}} \bigl[\mu(\*s, \*z) + \rho(\*s, \*x')\bigr]
\le
\inf_{\*s \in \3I} \bigl[\mu(\*s, \*z) + \rho(\*s, \*x')\bigr].
\]
Taking the minimum with $\overline{y}$ on both sides gives $U_{\*z} \le U_{\*z, \3I}$, and symmetrically $L_{\*z, \3I} \le L_{\*z}$. Validity of Eq.~\ref{eq:finite_bound} then follows from Thm.~\ref{thm:cate_bound} and monotonicity of the expectation:
\[
\3E\bigl[L_{\*z, \3I}(\*x, \*X') \mid \*z\bigr]
\le
\3E\bigl[L_{\*z}(\*x, \*X') \mid \*z\bigr]
\le
\tau(\*x, \*z)
\le
\3E\bigl[U_{\*z}(\*x, \*X') \mid \*z\bigr]
\le
\3E\bigl[U_{\*z, \3I}(\*x, \*X') \mid \*z\bigr].
\]
It also follows directly from Lem.~\ref{lem:rect_anchor} at the points of $\3I$, without Eq.~\ref{eq:consistency}. Enlarging $\3I$ enlarges the index sets of the infimum and supremum, so the same comparison shows that the bounds tighten as $\3I$ grows. Finally, when $P(\*X \mid \*z)$ is an empirical law on $n$ points, each of the $n$ received treatments requires a minimum and a maximum over the $|\3I|$ reference points, which gives the operation count $O(n|\3I|)$. This completes the proof.
\end{proof}

\begin{corollary}\label{cor:strict}
Fix $\*x \in \3X$ and a reference set $\3I \subseteq S_{\*z}$. Then $\3E[U_{\*z}(\*x, \*X') \mid \*z] < \3E[U_{\*z, \3I}(\*x, \*X') \mid \*z]$ if and only if $U_{\*z}(\*x, \*X') < U_{\*z, \3I}(\*x, \*X')$ with positive probability, that is, if and only if some reference point $\*s \in S_{\*z}$ improves on the minimum over $\3I$ (and on $\overline{y}$) for a set of received treatments of positive probability. The same holds for the lower bounds.
\end{corollary}
\begin{proof}
By Cor.~\ref{cor:finite}, the gap
\[
D(\*x') \triangleq U_{\*z, \3I}(\*x, \*x') - U_{\*z}(\*x, \*x') \ge 0
\]
is nonnegative, and it is bounded. The expectation of a nonnegative function is positive if and only if the function is positive with positive probability, so
\[
\3E\bigl[U_{\*z}(\*x, \*X') \mid \*z\bigr] < \3E\bigl[U_{\*z, \3I}(\*x, \*X') \mid \*z\bigr]
\iff
P\bigl(D(\*X') > 0 \bigm| \*z\bigr) > 0.
\]
By the definition Eq.~\ref{eq:envelopes}, $D(\*x') > 0$, that is, $U_{\*z}(\*x, \*x') < U_{\*z, \3I}(\*x, \*x')$, holds exactly when some $\*s \in S_{\*z}$ has
\[
\mu(\*s, \*z) + \rho(\*s, \*x')
<
\overline{y} \wedge \inf_{\*s'' \in \3I}\bigl[\mu(\*s'', \*z) + \rho(\*s'', \*x')\bigr].
\]
The lower bounds are symmetric. This completes the proof.
\end{proof}

\begin{corollary}\label{cor:manski}
Let $\3X$ be finite with $\3Y = [\underline{y}, \overline{y}]$, and write $\pi(\*x, \*z) \triangleq P(\*X = \*x \mid \*z)$. If $\rho_{\mathrm{T}} = \rho_{\mathrm{S}} = c\,\mathds{1}\{\cdot \neq \cdot\}$ with $c \ge \operatorname{diam}(\3Y)$ (no smoothness), then for every $\*x \in \3X$ the sharp set $\Theta(\*x, \*z)$ of Thm.~\ref{thm:sharp} (Eq.~\ref{eq:identified_set}) is the conditional Manski interval \citep{manski1990nonparametric}
\begin{equation}\label{eq:manski}
\begin{split}
    \Theta(\*x, \*z) = \bigl[\, &\pi(\*x, \*z)\,\mu(\*x, \*z) + \bigl(1 - \pi(\*x, \*z)\bigr) \underline{y}, \\
    &\pi(\*x, \*z)\,\mu(\*x, \*z) + \bigl(1 - \pi(\*x, \*z)\bigr) \overline{y} \,\bigr],
\end{split}
\end{equation}
with $\pi\mu$ read as $0$ when $\pi(\*x, \*z) = 0$. For $\*x \in S_{\*z}$ and any $\rho_{\mathrm{T}}$, the singleton reference set $\3I = \{\*x\}$ of Cor.~\ref{cor:finite} with $\rho_{\mathrm{S}} = c\,\mathds{1}\{\cdot \neq \cdot\}$ returns the same interval.
\end{corollary}
\begin{proof}
The costs $c\,\mathds{1}\{\cdot \neq \cdot\}$ are pseudometrics, and they are continuous on the finite set $\3X$. Moreover, for all $\*s \neq \*s'$,
\[
|\mu(\*s, \*z) - \mu(\*s', \*z)|
\le
\operatorname{diam}(\3Y)
\le
2c
=
\rho_{\mathrm{T}}(\*s, \*s') + \rho_{\mathrm{S}}(\*s, \*s').
\]
So Eq.~\ref{eq:consistency} holds, and Thm.~\ref{thm:sharp} applies.

We next compute the bounds. Consider the terms $\mu(\*s, \*z) + \rho_{\mathrm{T}}(\*x, \*s) + \rho_{\mathrm{S}}(\*x', \*s)$ of Eq.~\ref{eq:envelopes}. If $\*x' = \*x \in S_{\*z}$, the reference point $\*s = \*x$ gives the term $\mu(\*x, \*z)$, and Thm.~\ref{thm:sharp}(a) gives
\[
U_{\*z}(\*x, \*x) = \mu(\*x, \*z).
\]
If $\*x' \neq \*x$, every $\*s$ differs from $\*x$ or from $\*x'$, so at least one of the two costs equals $c$, and every term satisfies
\[
\mu(\*s, \*z) + \rho_{\mathrm{T}}(\*x, \*s) + \rho_{\mathrm{S}}(\*x', \*s)
\ge
\underline{y} + c
\ge
\overline{y}.
\]
Hence $U_{\*z}(\*x, \*x') = \overline{y}$ for $\*x' \neq \*x$. Symmetrically,
\[
L_{\*z}(\*x, \*x) = \mu(\*x, \*z),
\qquad
L_{\*z}(\*x, \*x') = \underline{y} \quad (\*x' \neq \*x).
\]

Finally, we average over $\*X' \sim P(\*X \mid \*z)$. The atom $\*X' = \*x$ has probability $\pi(\*x, \*z)$ (zero if $\*x \notin S_{\*z}$), and every other received treatment takes the trivial values. Therefore
\[
\3E\bigl[U_{\*z}(\*x, \*X') \mid \*z\bigr] = \pi(\*x, \*z)\,\mu(\*x, \*z) + \bigl(1 - \pi(\*x, \*z)\bigr)\overline{y},
\qquad
\3E\bigl[L_{\*z}(\*x, \*X') \mid \*z\bigr] = \pi(\*x, \*z)\,\mu(\*x, \*z) + \bigl(1 - \pi(\*x, \*z)\bigr)\underline{y},
\]
and Thm.~\ref{thm:sharp}(b) gives Eq.~\ref{eq:manski}.

For the singleton $\3I = \{\*x\}$, the only term is
\[
\mu(\*x, \*z) + \rho_{\mathrm{T}}(\*x, \*x) + \rho_{\mathrm{S}}(\*x', \*x)
=
\mu(\*x, \*z) + c\,\mathds{1}\{\*x' \neq \*x\},
\]
whatever $\rho_{\mathrm{T}}$ is. It equals $\mu(\*x, \*z)$ at $\*x' = \*x$ and is at least $\underline{y} + c \ge \overline{y}$ otherwise, so the same computation applies and returns the same interval. This completes the proof.
\end{proof}

\subsection{Sharpness analysis}\label{app:proofs_sharp}
Call a function $g : \3X \times S_{\*z} \to \3Y$ \emph{admissible} if it satisfies the inequalities of \ref{ass:str} and \ref{ass:sts} and has diagonal $\mu(\cdot, \*z)$. By Lem.~\ref{lem:product}, such a diagonal satisfies Eq.~\ref{eq:consistency} and is continuous, so we take the continuous version of $\mu(\cdot, \*z)$. By consistency and \ref{ass:bounded}, the ETT function of every SCM in $\2M(\rho_{\mathrm{T}}, \rho_{\mathrm{S}})$ that reproduces $P(\*X, \*Z, Y)$ is admissible. The following lemma establishes the converse.

\begin{lemma}[Realization]\label{lem:realization}
Let $\{g_{\*z}\}_{\*z \in \3Z}$ be a family of admissible functions, jointly measurable in $(\*z, \*x, \*x')$. Then there exists an SCM in $\2M(\rho_{\mathrm{T}}, \rho_{\mathrm{S}})$ that reproduces $P(\*X, \*Z, Y)$ and whose ETT function is $g_{\*z}$ at every covariate value $\*z$.
\end{lemma}
\begin{proof}
Let $(\tilde{\*Z}, \tilde{\*X}) \sim P(\*Z, \*X)$ and $\varepsilon \sim \mathrm{Unif}(0, 1)$ be independent, and let $Q_{\*x', \*z}$ be the quantile function of $Y$ given $(\*X, \*Z) = (\*x', \*z)$, chosen jointly measurable. Take the latent confounder $\*U \triangleq (\tilde{\*Z}, \tilde{\*X})$, the outcome noise $\varepsilon_Y \triangleq \varepsilon$, a degenerate $\varepsilon_X$, and set $\*Z \gets \tilde{\*Z}$, $\*X \gets \tilde{\*X}$, and
\begin{equation*}
    Y \gets f_Y(\*X, \*Z, \*U, \varepsilon_Y), \qquad f_Y(\*x, \*z, (\tilde{\*z}, \tilde{\*x}), \varepsilon) \triangleq
    \begin{cases}
        Q_{\tilde{\*x}, \tilde{\*z}}(\varepsilon), & \*x = \tilde{\*x}, \\
        g_{\tilde{\*z}}(\*x, \tilde{\*x}), & \*x \neq \tilde{\*x}.
    \end{cases}
\end{equation*}

We first verify that this SCM reproduces the observational law. Observationally, $\*X = \tilde{\*X}$, so the first branch of $f_Y$ applies, and
\[
(\*Z, \*X) = (\tilde{\*Z}, \tilde{\*X}) \sim P(\*Z, \*X),
\qquad
Y = Q_{\*X, \*Z}(\varepsilon).
\]
Since $\varepsilon \sim \mathrm{Unif}(0, 1)$ is independent of $(\tilde{\*Z}, \tilde{\*X})$, $Q_{\*X, \*Z}(\varepsilon)$ has the true conditional law of $Y$ given $(\*X, \*Z)$. Hence $P(\*X, \*Z, Y)$ is reproduced. Moreover, all values lie in $\3Y$: the conditional quantiles of $Y$ lie in $\3Y$ by \ref{ass:bounded}, and $g_{\tilde{\*z}}$ takes values in $\3Y$ by admissibility.

It remains to verify that the ETT function of this SCM is $g_{\*z}$. Under $\doo(\*x)$, the potential outcome is
\[
Y_{\*x} = f_Y(\*x, \tilde{\*Z}, \*U, \varepsilon),
\]
and conditioning on the received treatment $\*x'$ and the covariate value $\*z$ amounts to conditioning on $(\tilde{\*Z}, \tilde{\*X}) = (\*z, \*x')$. For $\*x \neq \*x'$, the second branch applies, and
\[
\3E[Y_{\*x} \mid \*x', \*z] = g_{\*z}(\*x, \*x').
\]
For $\*x = \*x'$, the first branch applies, and since $\varepsilon$ is independent of $(\tilde{\*Z}, \tilde{\*X})$,
\[
\3E[Y_{\*x'} \mid \*x', \*z] = \3E\bigl[Q_{\*x', \*z}(\varepsilon)\bigr] = \mu(\*x', \*z) = g_{\*z}(\*x', \*x'),
\]
where the last equality uses that $g_{\*z}$ has diagonal $\mu(\cdot, \*z)$. So the ETT function of this SCM is $g_{\*z}$ at every covariate value. Since each $g_{\*z}$ is admissible, it satisfies \ref{ass:str} and \ref{ass:sts}, and the SCM lies in $\2M(\rho_{\mathrm{T}}, \rho_{\mathrm{S}})$. This completes the proof.
\end{proof}

\thmsharp*
\begin{proof}
We first prove (a). For each $\*s \in S_{\*z}$, consider the map
\[
\varphi_{\*s}(\*p) \triangleq \mu(\*s, \*z) + d_{\*z}\bigl(\*p, (\*s, \*s)\bigr),
\qquad \*p \in \3X \times S_{\*z}.
\]
By the triangle inequality for $d_{\*z}$, each $\varphi_{\*s}$ is $d_{\*z}$-Lipschitz:
\[
|\varphi_{\*s}(\*p) - \varphi_{\*s}(\*q)|
=
\bigl|d_{\*z}\bigl(\*p, (\*s, \*s)\bigr) - d_{\*z}\bigl(\*q, (\*s, \*s)\bigr)\bigr|
\le
d_{\*z}(\*p, \*q).
\]
An infimum of $d_{\*z}$-Lipschitz functions is $d_{\*z}$-Lipschitz whenever it is finite, and here it is bounded below:
\[
\inf_{\*s \in S_{\*z}} \varphi_{\*s}(\*p) \ge \inf_{\*s \in S_{\*z}} \mu(\*s, \*z) \ge \underline{y}.
\]
Taking the minimum with the constant $\overline{y}$ preserves the property. So $U_{\*z} = \bigl(\inf_{\*s} \varphi_{\*s}\bigr) \wedge \overline{y}$ is $d_{\*z}$-Lipschitz and takes values in $\3Y$, and by Lem.~\ref{lem:product} it satisfies the inequalities of \ref{ass:str} and \ref{ass:sts}.

It remains to compute $U_{\*z}$ on the diagonal. Fix $\*s \in S_{\*z}$. The reference point $\*s$ itself gives
\[
U_{\*z}(\*s, \*s) \le \mu(\*s, \*z) + \rho_{\mathrm{T}}(\*s, \*s) + \rho_{\mathrm{S}}(\*s, \*s) = \mu(\*s, \*z).
\]
Conversely, by Eq.~\ref{eq:consistency}, for every $\*s' \in S_{\*z}$,
\[
\mu(\*s', \*z) + \rho_{\mathrm{T}}(\*s, \*s') + \rho_{\mathrm{S}}(\*s, \*s') \ge \mu(\*s, \*z),
\]
and $\overline{y} \ge \mu(\*s, \*z)$. Taking the infimum over $\*s'$ and the minimum with $\overline{y}$ gives $U_{\*z}(\*s, \*s) \ge \mu(\*s, \*z)$. Hence $U_{\*z}(\*s, \*s) = \mu(\*s, \*z)$. The argument for $L_{\*z}$ is symmetric. This proves (a).

We next prove (b). By (a), $U_{\*z}$ and $L_{\*z}$ are admissible. The admissibility constraints are convex in $g$, so for every $\lambda \in [0, 1]$,
\[
g_\lambda \triangleq \lambda U_{\*z} + (1 - \lambda) L_{\*z}
\]
is admissible. The integrals $\3E[g_\lambda(\*x, \*X') \mid \*z]$ exist, because $g_\lambda(\*x, \cdot)$ is $\rho_{\mathrm{S}}$-Lipschitz, hence continuous and bounded. By linearity,
\[
\3E\bigl[g_\lambda(\*x, \*X') \mid \*z\bigr]
=
\lambda\, \3E\bigl[U_{\*z}(\*x, \*X') \mid \*z\bigr] + (1 - \lambda)\, \3E\bigl[L_{\*z}(\*x, \*X') \mid \*z\bigr],
\]
which sweeps the interval of Eq.~\ref{eq:cate_two_sided} as $\lambda$ runs over $[0, 1]$. By Thm.~\ref{thm:cate_bound}, $\Theta(\*x, \*z)$ lies inside that interval. Conversely, fix $\lambda \in [0, 1]$ and take the family $g_\lambda$ built at every covariate value from the same $\lambda$. Since \ref{ass:str} and \ref{ass:sts} hold at every covariate value, $\mu(\cdot, \*z')$ satisfies Eq.~\ref{eq:consistency} for almost every $\*z'$, so by (a) this family is admissible at almost every covariate value. By Lem.~\ref{lem:realization}, it is the ETT function of an SCM $\1M \in \2M(\rho_{\mathrm{T}}, \rho_{\mathrm{S}})$ that reproduces $P(\*X, \*Z, Y)$, and by Eq.~\ref{eq:cate_expand},
\[
\3E_{\1M}[Y_{\*x} \mid \*z] = \3E\bigl[g_\lambda(\*x, \*X') \mid \*z\bigr].
\]
Hence every point of the interval lies in $\Theta(\*x, \*z)$, which proves (b).

Finally, we prove the refutation claim. If an SCM in $\2M(\rho_{\mathrm{T}}, \rho_{\mathrm{S}})$ reproduces the observational law, then for almost every $\*z$ its ETT function has diagonal $\mu(\cdot, \*z)$, which then satisfies Eq.~\ref{eq:consistency} by Lem.~\ref{lem:product}. Therefore, if Eq.~\ref{eq:consistency} fails on a set of covariate values of positive probability, no such SCM exists. This completes the proof.
\end{proof}

\begin{corollary}\label{cor:test}
Under \ref{ass:bounded}, the class $\2M(\rho_{\mathrm{T}}, \rho_{\mathrm{S}})$ contains an SCM that reproduces $P(\*X, \*Z, Y)$ if and only if Eq.~\ref{eq:consistency} holds for $P$-almost every $\*z$. With Euclidean costs this reads $|\mu(\*s, \*z) - \mu(\*s', \*z)| \le (L_Y + C)\|\*s - \*s'\|$ on $S_{\*z}$: the Lipschitz constant of the observed regression on the support is an identified lower bound on $L_Y + C$. In particular, under \ref{ass:bounded}--\ref{ass:smoothness_X}, Eq.~\ref{eq:consistency} with the costs of Prop.~\ref{prop:structural} is a testable implication.
\end{corollary}
\begin{proof}
We first prove necessity. If an SCM in $\2M(\rho_{\mathrm{T}}, \rho_{\mathrm{S}})$ reproduces $P(\*X, \*Z, Y)$, then Eq.~\ref{eq:consistency} holds for $P$-almost every $\*z$ by the refutation part of Thm.~\ref{thm:sharp}.

We next prove sufficiency. Suppose Eq.~\ref{eq:consistency} holds for $P$-almost every $\*z$. By Thm.~\ref{thm:sharp}(a), $U_{\*z}$ is admissible at almost every $\*z$, and Lem.~\ref{lem:realization}, applied to the family $\{U_{\*z}\}$, gives an SCM in $\2M(\rho_{\mathrm{T}}, \rho_{\mathrm{S}})$ that reproduces $P(\*X, \*Z, Y)$.

Finally, with Euclidean costs $\rho_{\mathrm{T}} = L_Y\|\cdot\|$ and $\rho_{\mathrm{S}} = C\|\cdot\|$,
\[
\rho_{\mathrm{T}}(\*s, \*s') + \rho_{\mathrm{S}}(\*s, \*s') = (L_Y + C)\|\*s - \*s'\|,
\]
so Eq.~\ref{eq:consistency} takes the stated form. This completes the proof.
\end{proof}

\propkim*
\begin{proof}
With $\3I = \{x, x'\}$ at the received treatment $x'$, Eq.~\ref{eq:envelopes} has two terms. The reference point $s = x$ contributes
\[
\mu(x) + b \cdot 0 + a|x' - x| = \mu(x) + a|x' - x|,
\]
and the reference point $s = x'$ contributes
\[
\mu(x') + b|x - x'| + a \cdot 0 = \mu(x') + b|x' - x|.
\]
The minimum of the two is the integrand of $u_2(x)$, so averaging over the received treatment gives the bound of Cor.~\ref{cor:finite} with this reference set. Validity follows from Lem.~\ref{lem:product} and consistency at these two points, pointwise in $x'$. This completes the proof.
\end{proof}

\propkimloose*
\begin{proof}
Take $a = b = 1$, the treatment uniform on $S = \{0, 1, 2\}$, $Y = \mu(X)$ with $\mu = (0, -1, 0)$, and the query $x = 2$. Write $g(x, x') \triangleq \3E[Y_x \mid X = x']$ for $x, x' \in S$. The data fix the diagonal,
\[
g(x', x') = \mu(x'), \qquad x' \in S.
\]
The two-point bound uses the columns $s \in \{2, x'\}$ of Table~\ref{tab:kim}, which give the row values $(2, 0, 0)$ at $x' = 0, 1, 2$. Hence
\[
u_2(2) = \tfrac{1}{3}(2 + 0 + 0) = \tfrac{2}{3}.
\]

We first show that every model has $\3E[Y_{x=2}] \le 1/3$. Take any model whose ETT satisfies STR and STS with $a = b = 1$ and whose diagonal is $\mu$. At $x' = 2$, the diagonal gives $g(2, 2) = \mu(2) = 0$. At $x' = 1$, STR gives
\[
g(2, 1) \le g(1, 1) + 1 = 0.
\]
At $x' = 0$, STR and then STS at $x = 1$ give
\[
g(2, 0) \le g(1, 0) + 1 \le g(1, 1) + 2 = 1.
\]
Therefore
\[
\3E[Y_{x=2}] = \tfrac{1}{3}\bigl(g(2, 0) + g(2, 1) + g(2, 2)\bigr) \le \tfrac{1}{3}(1 + 0 + 0) = \tfrac{1}{3}.
\]
This uses only the two Lipschitz conditions and the diagonal, not any bound on the outcomes, so it holds over the whole class of \citet{kim2018identification}.

We next show that the value $1/3$ is attained. Let $g$ be the function with rows indexed by $x$ and columns by $x'$,
\[
\bigl(g(x, x')\bigr)_{x, x' \in S}
=
\begin{pmatrix}
    0 & 0 & 1 \\
    0 & -1 & 0 \\
    1 & 0 & 0
\end{pmatrix};
\]
this is the upper bound $\inf_{s}[\mu(s) + |x - s| + |x' - s|]$. Adjacent entries of every row and every column differ by at most $1$, so $g$ satisfies STR and STS with $a = b = 1$, and its diagonal is $\mu$. The model $X \sim \mathrm{Unif}(S)$, $Y_x \triangleq g(x, X)$ has ETT $g$ and observational outcome $Y = g(X, X) = \mu(X)$, so it reproduces the observational law, and
\[
\3E[Y_{x=2}] = \tfrac{1}{3}(1 + 0 + 0) = \tfrac{1}{3}.
\]
Combining the two steps, the sharp upper bound is $1/3 < 2/3 = u_2(2)$: the two-point bound is loose, in the setting of \citet{kim2018identification} as well as in ours. This completes the proof.
\end{proof}

\thmcontrast*
\begin{proof}
Write $U_j \triangleq U_{\*z}(\*x_j, \cdot)$ and $L_j \triangleq L_{\*z}(\*x_j, \cdot)$ on $S_{\*z}$ for $j \in \{0, 1\}$, and $d \triangleq d_{\*z}$. By Thm.~\ref{thm:cate_bound} applied to the admissible function $U_{\*z}$ of Thm.~\ref{thm:sharp}(a), $L_j \le U_j$.

We first prove (a). Let $g$ be admissible and $\*x' \in S_{\*z}$. By Thm.~\ref{thm:cate_bound},
\[
g(\*x_1, \*x') - g(\*x_0, \*x') \le U_1(\*x') - L_0(\*x'),
\]
and by \ref{ass:str},
\[
g(\*x_1, \*x') - g(\*x_0, \*x') \le \rho_{\mathrm{T}}(\*x_1, \*x_0) = r.
\]
Hence $g(\*x_1, \*x') - g(\*x_0, \*x') \le \min\{U_1(\*x') - L_0(\*x'),\ r\}$. Integrating the smaller of the two over $P(\*X \mid \*z)$ and using Eq.~\ref{eq:cate_expand} at both treatment levels gives
\[
\tau(\*x_1, \*z) - \tau(\*x_0, \*z) \le \3E\bigl[\min\{U_1(\*X') - L_0(\*X'),\ r\} \bigm| \*z\bigr],
\]
the upper endpoint of Eq.~\ref{eq:contrast}. Exchanging $\*x_1$ and $\*x_0$ gives the lower one.

We next prove (b), starting with the upper endpoint. Define on $S_{\*z}$
\[
a \triangleq U_1,
\qquad
b \triangleq \max\{L_0,\ U_1 - r\},
\qquad\text{so that}\qquad
a - b = \min\{U_1 - L_0,\ r\}
\]
pointwise. It remains to exhibit an admissible $g$ with $g(\*x_1, \cdot) = a$ and $g(\*x_0, \cdot) = b$. Let
\[
A \triangleq (\{\*x_1, \*x_0\} \times S_{\*z}) \cup D,
\qquad
D \triangleq \{(\*s, \*s) : \*s \in S_{\*z}\},
\]
and define $g$ on $A$ by $a$ on $\{\*x_1\} \times S_{\*z}$, by $b$ on $\{\*x_0\} \times S_{\*z}$, and by $\mu(\cdot, \*z)$ on the diagonal $D$. The triangle inequality for $\rho_{\mathrm{T}}$ gives
\[
\rho_{\mathrm{T}}(\*x_1, \*s) \le r + \rho_{\mathrm{T}}(\*x_0, \*s),
\]
and the clipping at $\overline{y}$ and $\underline{y}$ preserves the resulting inequalities
\begin{equation}\label{eq:shift}
    U_1 \le U_0 + r, \qquad L_0 \le L_1 + r.
\end{equation}
We check that $g$ is $d$-Lipschitz on $A$, pair type by pair type.
\begin{itemize}[leftmargin=*,topsep=2pt,itemsep=1pt]
    \item \emph{Diagonal--diagonal}: this is Eq.~\ref{eq:consistency}.
    \item \emph{$a$--diagonal}: the requirement is
    \[
    |a(\*x') - \mu(\*s, \*z)| \le \rho_{\mathrm{T}}(\*x_1, \*s) + \rho_{\mathrm{S}}(\*x', \*s) \qquad \text{for all } \*s, \*x' \in S_{\*z},
    \]
    which holds whenever $L_1 \le a \le U_1$, by the definitions Eq.~\ref{eq:envelopes}; here $a = U_1$.
    \item \emph{$b$--diagonal}: likewise, this requires $L_0 \le b \le U_0$. The lower inequality holds by the definition of $b$, and the upper one because $L_0 \le U_0$ and $U_1 - r \le U_0$ by Eq.~\ref{eq:shift}.
    \item \emph{$a$--$a$ and $b$--$b$}: $U_1$, $L_0$, and $U_1 - r$ are $\rho_{\mathrm{S}}$-Lipschitz by Thm.~\ref{thm:sharp}(a), and a maximum of $\rho_{\mathrm{S}}$-Lipschitz functions is $\rho_{\mathrm{S}}$-Lipschitz.
    \item \emph{$a(\*s) - b(\*s') \le r + \rho_{\mathrm{S}}(\*s, \*s')$}: by the definition of $b$ and the $\rho_{\mathrm{S}}$-Lipschitz property of $U_1$,
    \[
    b(\*s') \ge U_1(\*s') - r \ge U_1(\*s) - \rho_{\mathrm{S}}(\*s, \*s') - r.
    \]
    \item \emph{$b(\*s') - a(\*s) \le r + \rho_{\mathrm{S}}(\*s, \*s')$}: if $b(\*s') = U_1(\*s') - r$, this follows from the $\rho_{\mathrm{S}}$-Lipschitz property of $U_1$. If $b(\*s') = L_0(\*s')$, then by Eq.~\ref{eq:shift},
    \[
    L_0(\*s') \le L_1(\*s') + r \le U_1(\*s') + r \le U_1(\*s) + \rho_{\mathrm{S}}(\*s, \*s') + r.
    \]
\end{itemize}
The values lie in $\3Y$, since $\underline{y} \le L_0 \le b \le U_0 \le \overline{y}$. Points of $A$ that coincide receive consistent values. If $\*x_1 \in S_{\*z}$, then $a(\*x_1) = U_1(\*x_1) = \mu(\*x_1, \*z)$ by Thm.~\ref{thm:sharp}(a). If $\*x_0 \in S_{\*z}$, then $L_0(\*x_0) = \mu(\*x_0, \*z)$, and the reference point $\*s = \*x_0$ in Eq.~\ref{eq:envelopes} gives
\[
U_1(\*x_0) - r \le \mu(\*x_0, \*z) + \rho_{\mathrm{T}}(\*x_1, \*x_0) + \rho_{\mathrm{S}}(\*x_0, \*x_0) - r = \mu(\*x_0, \*z),
\]
so $b(\*x_0) = \mu(\*x_0, \*z)$.

We now extend $g$ from $A$ to $\3X \times S_{\*z}$ by the McShane formula \citep{mcshane1934extension}
\[
\hat g(\*q) \triangleq \min\Bigl\{\overline{y},\ \inf_{\*p \in A}\bigl[g(\*p) + d(\*q, \*p)\bigr]\Bigr\},
\qquad \*q \in \3X \times S_{\*z}.
\]
It is $d$-Lipschitz, it agrees with $g$ on $A$ because $g$ is $d$-Lipschitz on $A$ and bounded by $\overline{y}$, and it takes values in $\3Y$ because $g \ge \underline{y}$ on $A$. By Lem.~\ref{lem:product}, $\hat g$ satisfies the inequalities of \ref{ass:str} and \ref{ass:sts}, and it has diagonal $\mu(\cdot, \*z)$ because $D \subseteq A$, so it is admissible. As in the proof of Thm.~\ref{thm:sharp}(b), $U_{\*z'}$ is admissible at almost every covariate value $\*z'$. By Lem.~\ref{lem:realization}, with $\hat g$ at covariate value $\*z$ and $U_{\*z'}$ at every other covariate value, there is an SCM in $\2M(\rho_{\mathrm{T}}, \rho_{\mathrm{S}})$ that reproduces the observational law, and by Eq.~\ref{eq:cate_expand} its CATE is
\[
\3E\bigl[a(\*X') - b(\*X') \bigm| \*z\bigr] = \3E\bigl[\min\{U_1(\*X') - L_0(\*X'),\ r\} \bigm| \*z\bigr],
\]
the upper endpoint.

Finally, exchanging the roles of $\*x_1$ and $\*x_0$ attains the lower endpoint. The admissible functions form a convex set and the CATE is linear in the ETT function, so the identified set is the whole interval. This completes the proof.
\end{proof}

%% file: appendix_b.tex
\section{Experimental Setups}\label{appendix:b}
This appendix specifies the three settings of Sec.~\ref{sec:experiments}: the V-shaped discrete SCMs (Sec.~\ref{sec:discrete_separation}), the continuous V-shaped and smooth V-shaped families (Sec.~\ref{sec:continuous_separation}), and the controlled IST-derived binary model (Sec.~\ref{sec:ist_controlled_variants}). It then reports two additional experiments: smooth polynomial, logistic, phase, and curved-linear separation examples (Sec.~\ref{sec:additional_continuous}), and a controlled reference-set ablation on the counterexample of Prop.~\ref{prop:kimloose} (Sec.~\ref{sec:kim_reference_ablation}). Conditional on each model, regressions and targets are finite latent sums; continuous bound integration uses the stated numerical quadrature. All experiments run on CPU on a single laptop (Apple M1 Pro, 16 GB memory); no GPUs are used.

\label{sec:experiment_protocol} Ours and Kim receive identical regressions, valid smoothness costs, outcome ranges, and integration weights. Kim uses the per-pair endpoints; ours adds the supported reference grid before clipping and averaging. Synthetic SCMs use their exact structural constants, with STR cost $L_Y|x-s|$ and STS cost $C|x'-s|$, where $C=L_X\operatorname{diam}(\mathcal Y)/p_{\min}$. Controlled IST instead uses exact finite-grid STR/STS constants of its constructed treatment-on-treated means. CMSM receives the same conditional outcome distribution.

CMSM uses the smallest single $\Lambda^\star\ge1$ whose interval contains every evaluated target within one SCM or IST age group. This target-based calibration gives full coverage by construction and is not an observable calibration rule. Let $w_{jq}=P(U=u_j\mid X=x_q)$ and $y_{jq}$ be the conditional outcome atoms; for the IST-derived model these are Bernoulli atoms with probabilities $1-\mu(x_q)$ and $\mu(x_q)$. For a target $t_q=\tau(x_q)$ define
\[
 A_q=\sum_j w_{jq}(y_{jq}-t_q)_+,
 \qquad B_q=\sum_j w_{jq}(t_q-y_{jq})_+.
\]
The minimum per-query tilt satisfies $\log\Lambda_q=\tfrac12|\log A_q-\log B_q|$. We evaluate these sums in log space and use one $\Lambda^\star=\max_q\Lambda_q$ per model, rather than choosing a separate reported parameter per query. Zero-mass atoms are removed; a constant conditional outcome equal to its target needs $\Lambda_q=1$. A nonconstant law cannot reach an extreme atom with a finite tilt, and such an input fails explicitly. The exporter records the threshold separately from the $10^{-10}$ slack added to $\log_{10}\Lambda^\star$ for numerical evaluation (no slack when $\Lambda^\star=1$). It checks coverage, failure just below a nontrivial threshold at its binding query, and that the threshold does not exceed the sufficient budget from both directions of $p(x\mid u)/p(x)$. Unit checks compare this inversion with forward numerical search and the Bernoulli formula. Coverage calibration uses the known target; it does not establish membership of the true assignment law in the calibrated CMSM.

We use a treatment-axis adaptation of the Lipschitz-envelope construction of \citet{khan2024offpolicy}. Its anchors are the exact observational means on the sorted union of the observation and supported query grids, and its constant is the maximum absolute adjacent secant slope on that union. Every query has an exact anchor, so the envelope is the singleton $\{\mu(x_q)\}$ independently of the constant. The implementation verifies this identity before removing floating-point roundoff. Coverage is evaluated against $\tau(x_q)$ and may fail because these models are confounded; we do not assert full coverage for this baseline. This is an assumption-mismatch comparison, not a counterexample to Khan's unconfoundedness-based guarantee. The resulting zero-width intervals cover $14.11\%$ of discrete targets and $4\%$ in each continuous family and IST group.

Figure~\ref{fig:section5_intervals} uses saved system $0$ from the discrete and smooth V-shaped continuous families and the middle IST age tertile, without selecting systems for improvement. Discrete treatments use interval bars and target markers; the remaining panels connect evaluated grid points. The table pools both continuous families even though the figure displays only the smooth family. Run \texttt{plot\_section5.py} to reproduce the figure from saved intervals without rerunning experiments.

\subsection{Discrete SCMs with a V-shaped response}\label{sec:discrete_separation}
The discrete models use a response with a pronounced interior minimum, so that global references can tighten Kim's per-pair endpoints. Let $U\in\{-1,+1\}$ with $P(U=-1)=p_-$, and let $X$ have an odd, equally spaced support of size $K$ on $[-1,1]$. Define
\[
 P(X=x_j\mid U=u)=\frac{1+\rho u x_j}{K},\qquad
 Y_x=\frac{2(b|x|+cu+c)}{b+2c}-1.
\]
The canonical coefficients are $K=7$, $b=1$, $c=1/2$, $\rho=0.1$, and $p_-=1/2$. Assignment is generated by an independent uniform variable through the conditional inverse CDF. Both latent types have positive probability, assignment has full overlap, and the outcome attains both $-1$ and $+1$ at supported treatments. Thus the supplied clipping range is exact, rather than enlarged to create an apparent improvement.

Both methods receive the maximum finite-support difference quotients of the response and assignment probabilities. Independently, these have the analytic values
\[
 L_Y^{\star}=\frac{2b}{b+2c},\qquad
 L_X^{\star}=\frac{\rho}{K},\qquad
 p_{\min}=\frac{1-\rho}{K},\qquad
 C=\frac{2\rho}{1-\rho}.
\]
We use these structural constants rather than substituting smaller direct ETT constants. At the canonical coefficients, $\tau(x)=|x|-1/2$ and $\mu(x)=|x|-1/2+0.05x$, so confounding is nonzero, with maximum gap $2.5\%$ of the outcome range. The global reference $s=0$ supplies the response minimum that Kim's per-pair endpoint references may omit. Adding this one reference to the endpoints reduces the normalized mean width from $0.13908$ to $0.11187$, a $19.57\%$ reduction, and gives the same interval as using every supported reference. The improvement is strict at six of the seven queries. An independent finite LP agrees at every query; this verifies sharpness in the induced STR/STS class, not sharpness under the stronger structural assumptions.

\input{table_discrete_separation}
\begin{figure*}[t]
\centering
\includegraphics[width=\linewidth]{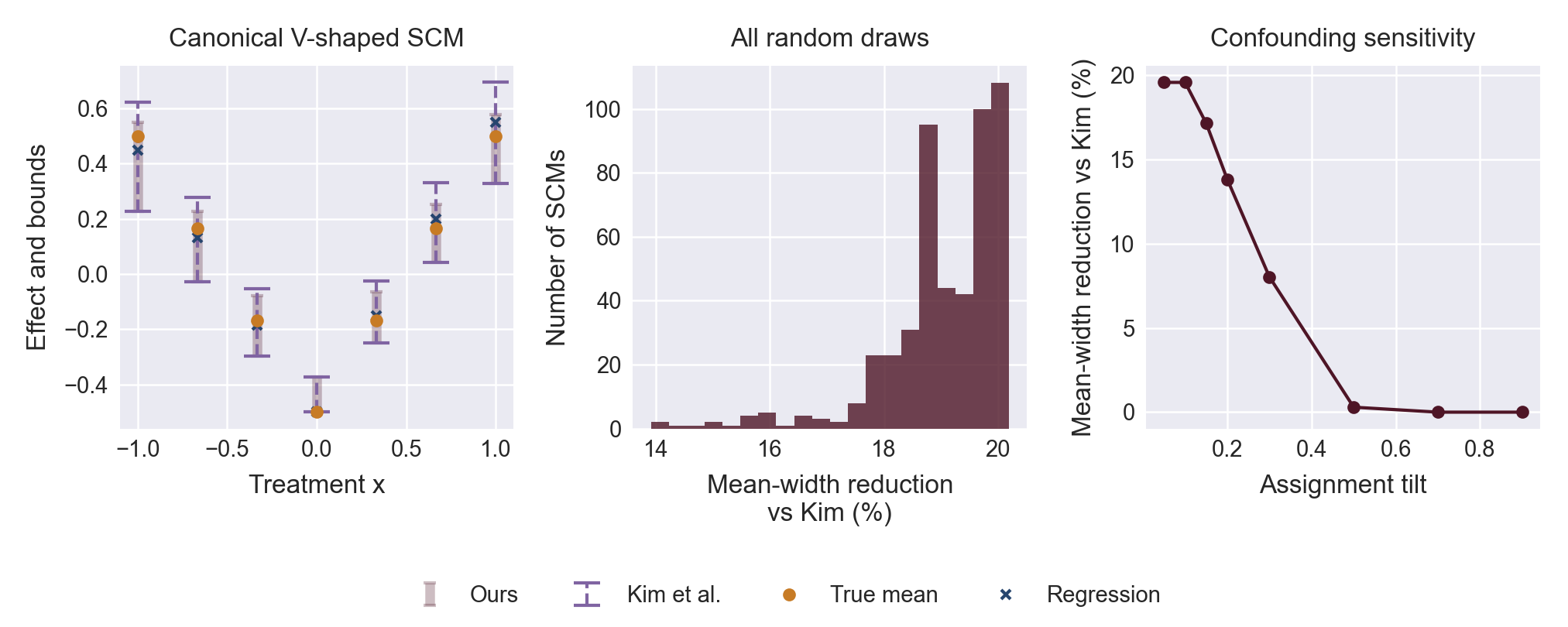}
\caption{Discrete V-shaped SCMs with matched structural constants and exact outcome-range clipping. Left: the canonical V-shaped response, its confounded observational regression, and both intervals. Middle: mean-width reductions for all $500$ random systems. Right: the reduction across the assignment-tilt sweep, including settings where the methods tie.}
\label{fig:discrete_separation}
\end{figure*}

We draw $500$ systems with $K$ uniform on $\{5,7,9\}$, $b\sim\texttt{Uniform}(0.9,1.1)$, $c\sim\texttt{Uniform}(0.45,0.55)$, $\rho\sim\texttt{Uniform}(0.08,0.12)$, and $p_-\sim\texttt{Uniform}(0.35,0.65)$, independently. This prior was designed after a feasibility check, so we do not interpret it as representative. No draw is rejected based on separation, coverage, or confounding. The pooled normalized widths are $0.11181$ for ours and $0.13832$ for Kim: a $19.17\%$ reduction over $3544$ queries, with per-system reductions ranging from $13.91\%$ to $20.20\%$ and full coverage for both methods (Table~\ref{tab:discrete_separation}).

A separate canonical tilt sweep uses $\rho\in\{0.05,0.1,0.15,0.2,0.3,0.5,0.7,0.9\}$, recalibrating both methods identically at every setting. At $\rho=0.2$, the confounding gap is $5\%$ of the outcome range and the mean-width reduction remains $13.82\%$. The improvement falls to zero at $\rho=0.7$ and $0.9$; coverage remains complete throughout (Fig.~\ref{fig:discrete_separation}). The script \texttt{lpi\_discrete\_separation.py} exports every random system, query interval, tilt setting, and nested reference ablation, together with the analytic-constant checks and canonical LP endpoints.

\begin{figure*}[t]
\centering
\includegraphics[width=\linewidth]{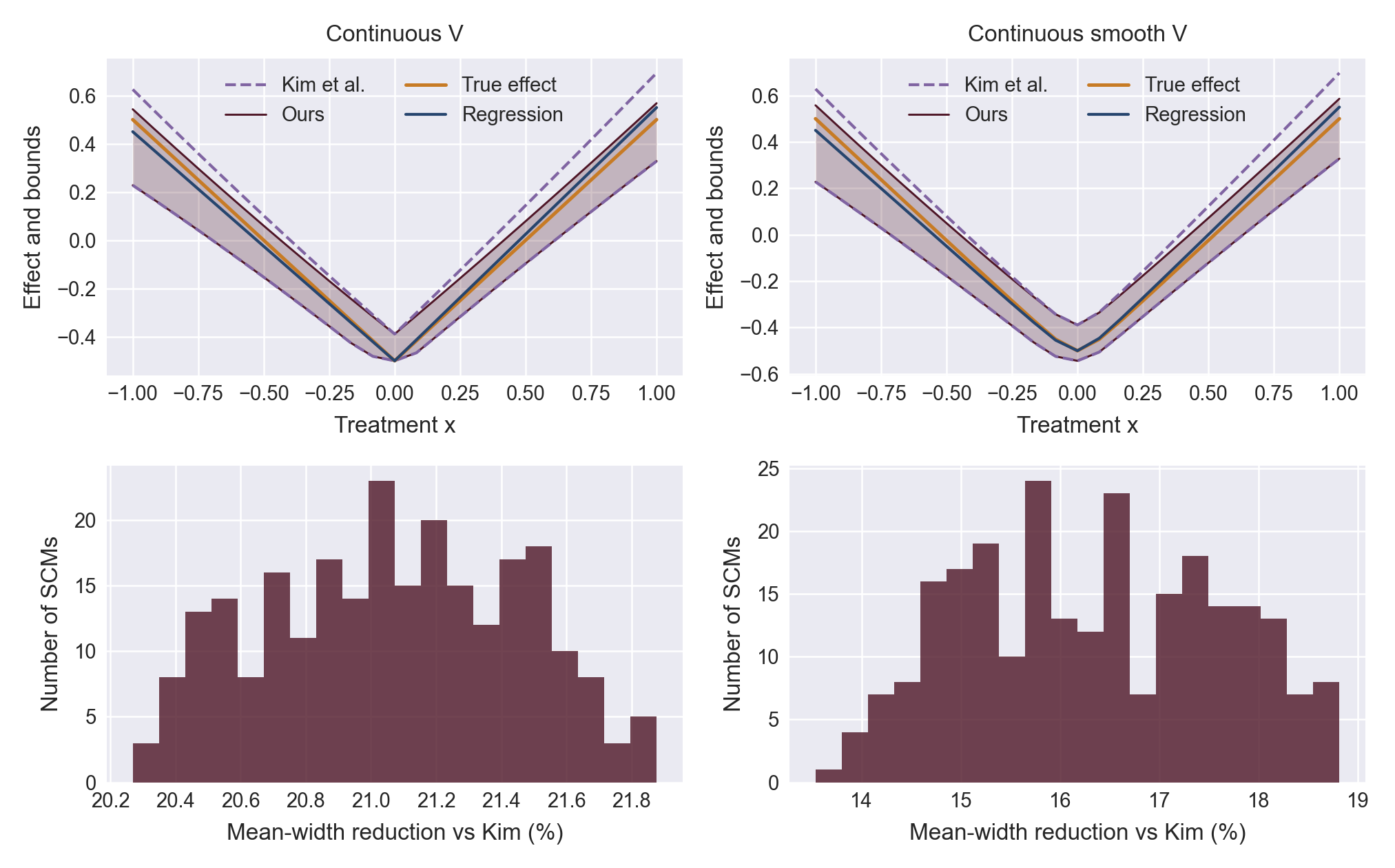}
\caption{Continuous V-shaped and smooth V-shaped response families. Top: canonical bounds at identical exact structural constants. Bottom: mean-width reductions against Kim across all $250$ random draws per family, including every draw without rejection.}
\label{fig:continuous_separation}
\end{figure*}

\subsection{Continuous SCMs with V-shaped responses}\label{sec:continuous_separation}
The continuous models extend the interior-minimum construction while retaining a strictly positive treatment density on the entire interval. Let $U\in\{-1,+1\}$, $P(U=-1)=p_-$, and define
\[
 p(x\mid u)=\frac{1+\rho ux}{2},\quad x\in[-1,1],
 \qquad Y=\frac{2\{b\,g_\epsilon(X)+cU+c\}}{b+2c}-1,
\]
where, writing $D=1+|m|$ and $N_\epsilon=\sqrt{D^2+\epsilon^2}-\epsilon$,
\[
 g_\epsilon(x)=\frac{\sqrt{(x-m)^2+\epsilon^2}-\epsilon}{N_\epsilon}.
\]
An independent uniform noise variable passed through the inverse conditional CDF defines the treatment equation. The response has its minimum at $m$ and maximum at the farther support endpoint, so $g_\epsilon$ attains $[0,1]$ and the outcome attains $[-1,1]$. With $b,c>0$, $0<\rho<1$, and $0<p_-<1$, the exact structural constants are
\[
 L_X=\rho/2,\qquad p_{\min}=(1-\rho)/2,\qquad
 L_Y=\frac{2bD}{(b+2c)N_\epsilon\sqrt{D^2+\epsilon^2}},
 \qquad C=\frac{2\rho}{1-\rho}.
\]
The case $\epsilon=0$ is a Lipschitz V-shaped response; $\epsilon>0$ gives a differentiable smooth approximation. Both methods use these same structural constants, with no substitution of smaller direct ETT constants.

The canonical systems use $b=1$, $c=0.5$, $\rho=0.1$, $p_-=0.5$, $m=0$, with $\epsilon=0$ or $0.05$. Independently draw $b\sim U[0.9,1.1]$, $c\sim U[0.45,0.55]$, $\rho\sim U[0.08,0.12]$, $p_-\sim U[0.35,0.65]$, and $m\sim U[-0.1,0.1]$. Set $\epsilon=0$ for the V-shaped family and draw $\epsilon\sim U[0.025,0.075]$ for the smooth family. Seed $0$ generates $250$ models per family, with every draw retained. This prior was selected after exploratory feasibility checks. Each model has $25$ equally spaced queries on $[-1,1]$. We use $257$ trapezoidal integration nodes and a $257$-point reference grid augmented by $m$ and the exact pairwise endpoint references. Regression and interventional means are exact sums over the two latent types.

Canonical mean-width reductions are $21.27\%$ and $16.27\%$; random-population reductions are $21.04\%$ and $16.34\%$, with full coverage for ours and Kim (Fig.~\ref{fig:continuous_separation}). Per-system reductions span $20.27$--$21.88\%$ and $13.54$--$18.81\%$; there is no acceptance threshold based on the improvement. All random V-shaped queries and $93.824\%$ of smooth-family queries improve strictly. Each model is checked for density normalization, overlap, attained outcome range, structural constants, coverage, and containment within Kim. Removing all additional references reproduces Kim at every query of all $502$ canonical and random systems. Canonical ablations add the minimum alone, then nested grids of $3,9,33,129,257,513,1025$ references augmented by the minimum, verifying interval nesting. A separate assignment-strength sweep retains $\rho\in\{0.05,0.1,0.15,0.2,0.3,0.5,0.7,0.9\}$, including weaker separation at stronger confounding.

For each family, refinement checks use the canonical system and the first and last random draws. We refine integration and reference grids jointly and separately through $257,513,1025$ points. Relative to the finest joint grid, the largest endpoint change is $1.674\times10^{-5}$ of the outcome range, below the declared $10^{-4}$ tolerance. Run \texttt{lpi\_continuous\_separation.py --emit-paper-assets --assert} to export parameters, intervals, summaries, reference ablations, assignment sweeps, refinement records, and the figure. The consolidated exporter then adds oracle CMSM and Khan to these same saved models.

\input{table_ist_controlled_variants}
\begin{figure*}[t]
\centering
\includegraphics[width=\linewidth]{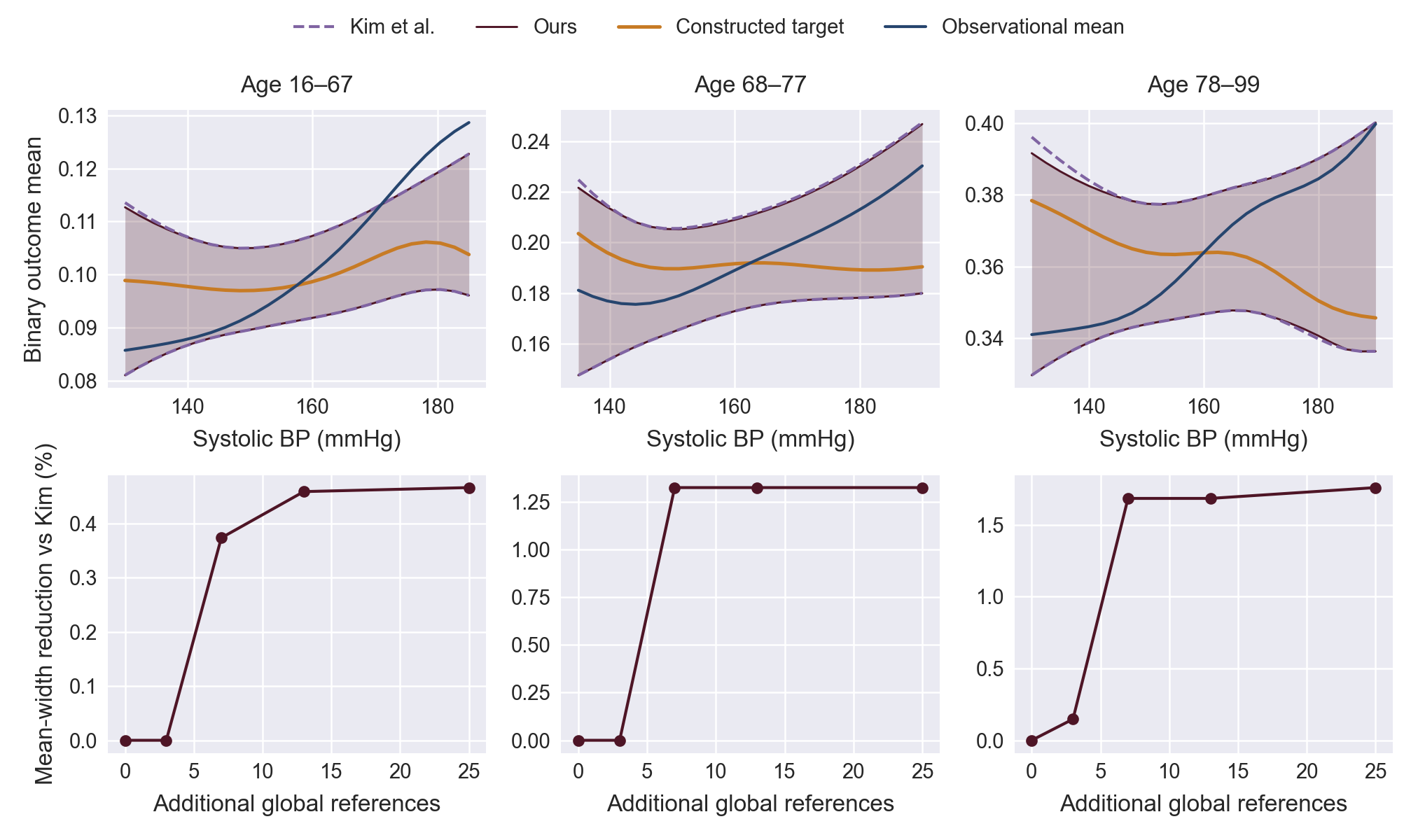}
\caption{Controlled IST-derived binary SCM, with identical oracle STR/STS constants for both methods. Top: bounds, constructed interventional targets, and observational means in each age group. Bottom: the reduction in mean width as global references are added to the per-pair endpoints.}
\label{fig:ist_controlled_variants}
\end{figure*}

\subsection{Controlled IST}\label{sec:ist_controlled_variants}
The International Stroke Trial \citep{carolei1997international}, whose anonymised patient-level database is publicly released on Edinburgh DataShare under the repository's end-user licence \citep{sandercock2011ist}, is used to construct a model with a known target, not to estimate clinical effects. Blood pressure is not randomized in IST, so its raw regression is not a causal dose--response. We take systolic blood pressure at randomization ($\texttt{RSBP}$) as the treatment, death within six months ($\texttt{FDEAD}$) as the binary outcome, the age tertile $z$ (roughly $6{,}400$ patients each) as the covariate, and consciousness at randomization $u\in\{\texttt{F},\texttt{D},\texttt{U}\}$ (fully alert $14{,}921$ / drowsy $4{,}254$ / unconscious $260$) as the hidden confounder. Consciousness strongly affects the outcome: the six-month death rates across its three levels are $0.150$, $0.461$, and $0.776$. Within each tertile, the empirical category frequencies define $P(U=u\mid z)$. There are $25$ equally spaced treatment values across the central $70\%$ of that tertile's observed systolic blood-pressure range, standardized to $t_j\in[-1,1]$. For each consciousness category, Nadaraya--Watson regression of the binary outcomes supplies probabilities $r_z(t_j,u)$. All categories within a tertile use one bandwidth, twice the full-tertile Silverman bandwidth in mmHg. These fitted probabilities define the new model; they are not assumed to identify clinical causal effects.

With $v(\texttt{F})=-1$ and $v(\texttt{D})=v(\texttt{U})=1$, set
\[
 P(X=t_j\mid U=u,z)=\frac{1+\eta v(u)t_j}{25},\qquad \eta=0.3.
\]
The symmetric grid makes this a normalized positive assignment law. An independent uniform variable transformed by the conditional inverse CDF defines $X$. A second, independent $E\sim\texttt{Uniform}(0,1)$ defines the binary potential outcome $Y_{t_j}=\mathbf{1}\{E\leq r_z(t_j,U)\}$, with target $\tau_z(t_j)=\sum_u P(u\mid z)r_z(t_j,u)$. The full treatment-on-treated matrix is known by Bayes' rule:
\[
 g_z(t_i,t_j)=\sum_u r_z(t_i,u)
 \frac{P(u\mid z)P(t_j\mid u,z)}{\sum_v P(v\mid z)P(t_j\mid v,z)},
 \qquad \mu_z(t_j)=g_z(t_j,t_j).
\]
Every target, regression, and integration weight is an exact finite sum conditional on the fitted model. No simulated outcomes or estimated effect targets enter the evaluation.

For each age group, compute the maximum difference quotient of $g_z$ over its first coordinate to obtain the STR coefficient $b_z$, and over its second coordinate to obtain the STS coefficient $a_z$. Both inequalities are checked over every matrix entry. These are oracle constants of the constructed mean-response model, not claims of pathwise smoothness for the binary outcome or of known propensity constants in the original data. Both methods receive exactly the same $a_z,b_z$, observational diagonal, treatment weights, and outcome range $[0,1]$. Our full reference set includes every supported treatment; Kim uses the two endpoint references at each pair $(t_i,t_j)$. Constants are passed through a direct-cost interface without inventing overlap or density parameters.

Table~\ref{tab:ist_controlled_variants} gives mean-width reductions of $0.47\%$, $1.32\%$, and $1.76\%$. Both methods cover all $75$ constructed means, and our intervals are contained in Kim's at every query. The constants and headline configuration were informed by a feasibility check, not preregistered. A nested reference ablation adds $0,3,7,13,25$ global points to the per-pair endpoints, reproducing Kim exactly at zero and tightening monotonically thereafter (Fig.~\ref{fig:ist_controlled_variants}). An independent sparse linear program over the $25^2$ ETT values, constrained only by the observed diagonal, STR, STS, and range, agrees with our full-reference endpoints to within $10^{-8}$ at the first, middle, and last query of every age group.

We retain every setting in sensitivity sweeps over $\eta\in\{0,0.1,0.3,0.5,0.7\}$; independent multipliers $\{1,1.25,2,4\}$ of the oracle STR and STS constants; bandwidth multipliers $\{1,2,3\}$ of the full-tertile Silverman bandwidth; and grid sizes $\{13,25,49\}$. The bandwidth and grid sweeps rebuild the model and its oracle constants. Both methods retain full coverage in all settings, including settings where the width improvement vanishes. At $\eta=0$, assignment is independent of the confounder, the STS coefficient is zero, and both bounds recover the model target exactly; relative width reductions are left undefined when Kim's width is zero. The script \texttt{lpi\_ist\_controlled\_variants.py} exports all settings, per-query intervals, fitted risks, latent proportions, calibration constants, and LP comparisons, with a checksum of the input data for reproducibility.

\input{table_additional_continuous_separation}
\begin{figure*}[t]
\centering
\includegraphics[width=\linewidth]{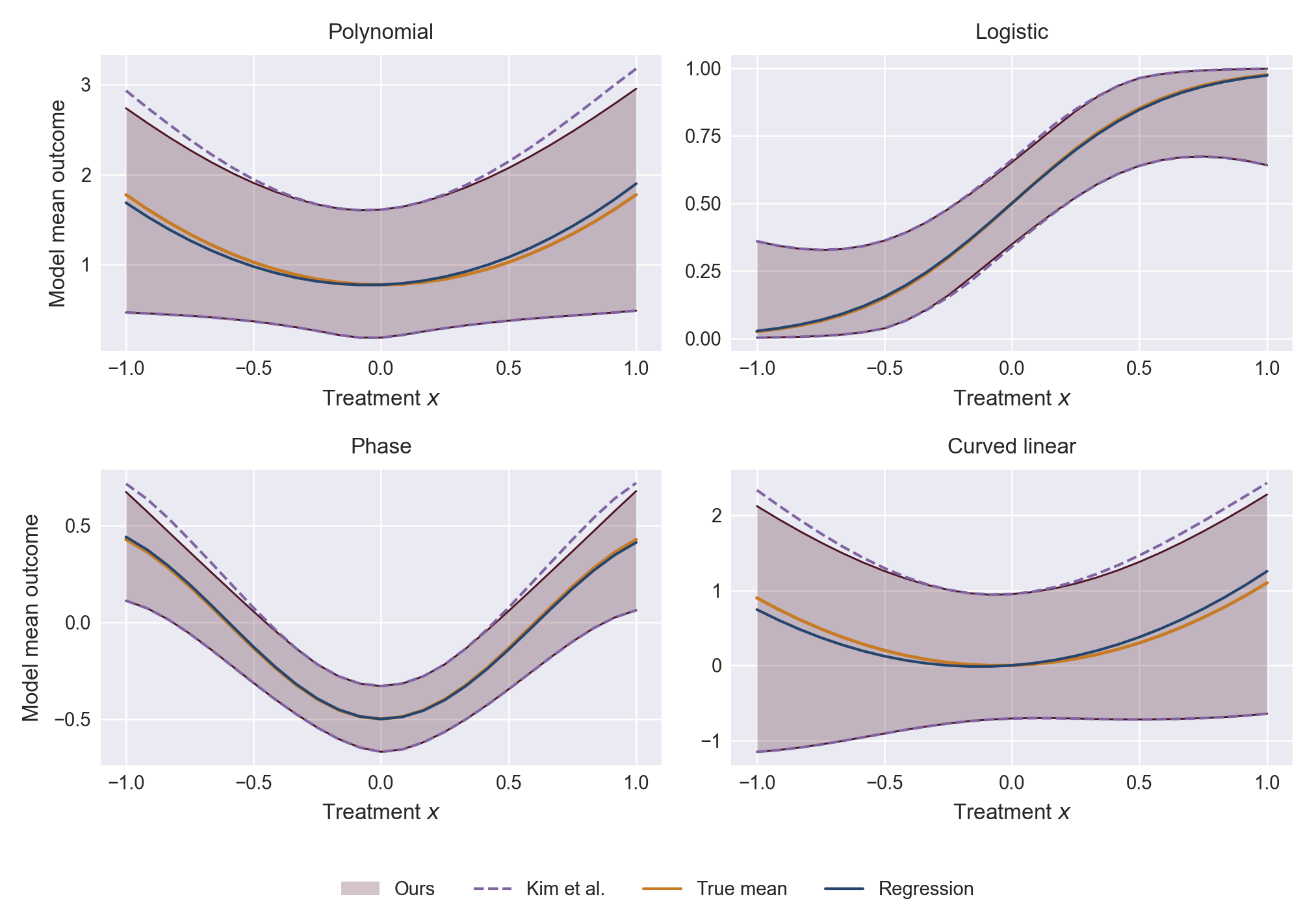}
\caption{Designed smooth continuous SCMs with identical structural constants and outcome ranges for both methods. Each panel shows our full-reference bounds (dark red), Kim's bounds (purple dashed), the true mean (orange), and the observational regression (blue). Top: polynomial and logistic; bottom: phase and curved linear. These are fixed canonical examples chosen after exploratory checks; Table~\ref{tab:additional_continuous_separation} also reports all random draws in the stated neighborhoods.}
\label{fig:additional_continuous_separation}
\end{figure*}

\subsection{Additional continuous SCMs}\label{sec:additional_continuous}

We reuse four smooth, non-V-shaped response families with a bounded-propensity assignment mechanism to demonstrate separation from Kim beyond the V-shaped examples. These canonical settings and their random neighborhoods were chosen after exploratory checks; they are designed demonstrations, not a representative prior or a claim of universal gains. Let $U$ have $201$ equally spaced atoms on $[-2,2]$, with probabilities proportional to $\exp(-u^2/2)$. On $[-1,1]$, set
\begin{align}
 p(x\mid u)&=\frac{1+\rho v(u)x}{2},\nonumber\\
 \text{Polynomial:}\quad Y_x&=b x^2+cU^2,
 & (b,c,\rho)&=(1,1,0.25),\nonumber\\
 \text{Logistic:}\quad Y_x&=\left(1+e^{-(bx-cU)}\right)^{-1},
 & (b,c,\rho)&=(4,1,0.25),\nonumber\\
 \text{Phase:}\quad Y_x&=\sin^2(wx)+\cos^2(wU)-1,
 & (w,\rho)&=(1.3,0.15),\nonumber\\
 \text{Curved linear:}\quad Y_x&=bx+qx^2+cU,
 & (b,q,c,\rho)&=(0.1,1,4,0.10).\nonumber
\end{align}
Here $v(u)=u^2/2-1$ for polynomial and phase, and $v(u)=u/2$ for logistic and curved linear. The curved-linear response includes a quadratic term and is not a strictly linear model.
Applying the conditional inverse CDF to an independent uniform variable defines $X$. Every conditional density is normalized and positive. Both methods receive the same exact structural constants $L_X=\rho/2$, $p_{\min}=(1-\rho)/2$, and $L_Y=2b$, $b/4$, $w$, or $|b|+2q$ for polynomial, logistic, phase, or curved linear, respectively, together with the exact outcome range over $x\in[-1,1]$ and the finite latent support. Thus the STS coefficient is $C=L_X\operatorname{diam}(\3Y)/p_{\min}$. Regressions and interventional means are exact latent sums. Both bounds use the same $257$-point trapezoidal integration grid and $25$ equally spaced queries; ours additionally uses all $257$ reference points, augmented by the per-pair endpoints.

With seed $0$, processing polynomial, phase, logistic, then curved linear, we draw $125$ systems per family. Each response coefficient is multiplied independently by $\texttt{Uniform}(0.9,1.1)$: $(b,c)$ for polynomial and logistic, $w$ for phase, and $(b,q,c)$ for curved linear. Independently, $\rho$ is drawn from $[0.20,0.30]$ for polynomial and logistic, $[0.10,0.20]$ for phase, or $[0.05,0.15]$ for curved linear. All $500$ draws are retained, without filtering by coverage, confounding, or separation. Table~\ref{tab:additional_continuous_separation} reports mean interval widths divided by each system's outcome range, with reductions computed from pooled mean widths. In the displayed order (polynomial, logistic, phase, curved linear), the canonical reductions are $4.37\%$, $1.18\%$, $6.29\%$, and $3.36\%$; the random-family reductions are $4.01\%$, $1.08\%$, $5.97\%$, and $2.51\%$. Both methods cover every query. Figure~\ref{fig:additional_continuous_separation} shows the bounds for the four fixed canonical examples.

Checks verify conditional-density normalization, structural and ETT smoothness costs, coverage, interval ordering, and containment within Kim. Endpoint-only references reproduce Kim; adding nested reference grids tightens monotonically. For each canonical example and the first and last random draw in each family, integration and reference grids are refined jointly and separately through $257$, $513$, and $1025$ points, with endpoint changes required to stay below $10^{-4}$ of the outcome range. Run \texttt{lpi\_additional\_continuous.py --emit-paper-assets} to reproduce the parameters, intervals, aggregate results, reference ablations, refinement records, table, and figure; \texttt{--plot-only} redraws the saved results without rerunning the experiment.

\subsection{Controlled reference-set ablation}\label{sec:kim_reference_ablation}
To isolate the role of additional reference points, we evaluate the counterexample of Prop.~\ref{prop:kimloose} with the observational law and both costs held fixed: $X$ is uniform on $\{0,1,2\}$, $\mu=(0,-1,0)$, and $a=b=1$. The supplied outcome range is $[-2,2]$, shared by both methods; clipping to this range leaves the two-point envelopes unchanged. This experiment concerns the STR+STS class directly, rather than the structural-overlap constants of the synthetic SCMs above. We use the attaining SCM from App.~\ref{appendix:a}, with $U$ uniform, $X=U$, and $Y_x=g(x,U)$, whose ETT rows are $(0,0,1)$, $(0,-1,0)$, and $(1,0,0)$. Its true interventional means are $(1/3,-1/3,1/3)$, and its STR and STS constants are both exactly one.

\begin{figure}[t]
\centering
\includegraphics[width=0.9\linewidth]{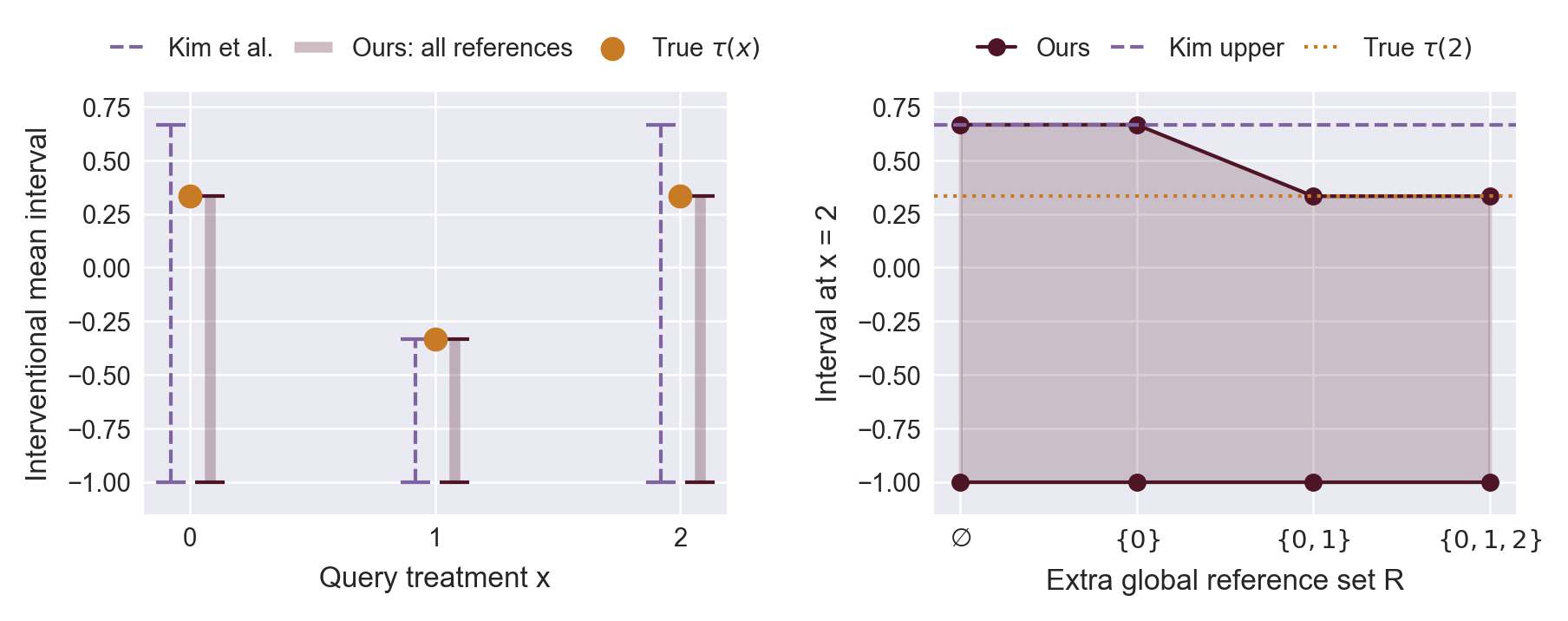}
\caption{Controlled reference ablation with identical data and constants. Left: Kim's intervals and our full-reference intervals at all three queries, with the true means marked. Right: the interval at $x=2$ as the extra global reference set $R$ grows. The interior reference $s=1$ lowers the upper endpoint to the attained sharp value $1/3$; adding the remaining endpoint does not tighten it further.}
\label{fig:kim_reference_ablation}
\end{figure}

At each pair $(x,x')$, our reference set is $\{x,x'\}\cup R$, with $R$ increasing through $\varnothing$, $\{0\}$, $\{0,1\}$, and $\{0,1,2\}$. With $R=\varnothing$, our implementation reproduces Kim exactly. At $x=2$, adding $s=0$ does not help; adding the interior reference $s=1$ lowers the upper endpoint from $2/3$ to $1/3$ while leaving the lower endpoint at $-1$. Width therefore falls from $5/3$ to $4/3$, a $20\%$ reduction, and all three true means remain covered at every stage (Table~\ref{tab:kim_reference_ablation}, Fig.~\ref{fig:kim_reference_ablation}). The final intervals agree with an independent linear program over all nine ETT values subject to the observed diagonal, STR, STS, and the outcome range. This is a controlled demonstration of strict improvement, not a claim of strict improvement on every system. A narrower supplied outcome range adds information and can remove the gap.

\input{table_kim_reference_ablation}

%% file: table_discrete_separation.tex
\begin{table}[t]
\centering
\caption{Stronger discrete variant with a V-shaped response. Widths are divided by the attained outcome range, $\operatorname{diam}(\mathcal{Y})=2$. Ours and Kim use identical exact structural constants; all random draws are retained.}
\label{tab:discrete_separation}
\begin{tabular}{lrrrrrr}
\toprule
Setting & SCMs & Queries & Ours width & Kim width & Reduction (\%) & Coverage (both) \\
\midrule
Canonical & 1 & 7 & 0.1119 & 0.1391 & 19.57 & 1.00/1.00 \\
Random & 500 & 3544 & 0.1118 & 0.1383 & 19.17 & 1.00/1.00 \\
\bottomrule
\end{tabular}
\end{table}

%% file: table_ist_controlled_variants.tex
\begin{table}[t]
\centering
\caption{Controlled IST-derived binary SCM: mean widths and coverage at identical oracle STR/STS constants. Each age group has 25 treatment queries. This model has a new, explicitly constructed target.}
\label{tab:ist_controlled_variants}
\begin{tabular}{lrrrrr}
\toprule
Age & Ours width & Kim width & Reduction (\%) & Ours cov. & Kim cov. \\
\midrule
16--67 & 0.01972 & 0.01981 & 0.47 & 1.00 & 1.00 \\
68--77 & 0.04844 & 0.04909 & 1.32 & 1.00 & 1.00 \\
78--99 & 0.04362 & 0.04440 & 1.76 & 1.00 & 1.00 \\
\bottomrule
\end{tabular}
\end{table}

%% file: table_additional_continuous_separation.tex
\begin{table}[t]
\centering
\caption{Designed smooth continuous SCMs. Widths are normalized by each system's outcome range. Random rows retain all 125 draws per family; reduction compares pooled mean widths. Canonical settings and neighborhoods were chosen after exploratory checks.}
\label{tab:additional_continuous_separation}
\begin{tabular}{llrrrrr}
\toprule
Setting & Family & Ours & Kim & Ours cov. & Kim cov. & Reduction (\%) \\
\midrule
Canonical & Polynomial & 0.3480 & 0.3639 & 1.00 & 1.00 & 4.37 \\
Canonical & Logistic & 0.3243 & 0.3282 & 1.00 & 1.00 & 1.18 \\
Canonical & Phase & 0.2136 & 0.2279 & 1.00 & 1.00 & 6.29 \\
Canonical & Curved linear & 0.1316 & 0.1362 & 1.00 & 1.00 & 3.36 \\
Random & Polynomial & 0.3471 & 0.3616 & 1.00 & 1.00 & 4.01 \\
Random & Logistic & 0.3235 & 0.3271 & 1.00 & 1.00 & 1.08 \\
Random & Phase & 0.2172 & 0.2310 & 1.00 & 1.00 & 5.97 \\
Random & Curved linear & 0.1244 & 0.1276 & 1.00 & 1.00 & 2.51 \\
\bottomrule
\end{tabular}
\end{table}

%% file: table_kim_reference_ablation.tex
\begin{table}[t]
\centering
\caption{Controlled reference ablation: fixed observational law, $a=b=1$, and known outcome range $[-2,2]$. Endpoints are $\{x,x'\}$ per stratum. Bounds and widths are at $x=2$ in outcome units; coverage is over all three queries. Adding $s=1$ reduces width by $20\%$ without losing coverage.}
\label{tab:kim_reference_ablation}
\begin{tabular}{lrrrr}
\toprule
References / method & Lower & Upper & Width & Coverage \\
\midrule
Kim et al. & -1.000 & 0.667 & 1.667 & 1.00 \\
Ours: endpoints & -1.000 & 0.667 & 1.667 & 1.00 \\
Ours: endpoints + $\{0\}$ & -1.000 & 0.667 & 1.667 & 1.00 \\
Ours: endpoints + $\{0,1\}$ & -1.000 & 0.333 & 1.333 & 1.00 \\
Ours: all references & -1.000 & 0.333 & 1.333 & 1.00 \\
\bottomrule
\end{tabular}
\end{table}